\documentclass[11pt]{article}

\usepackage{color}
\usepackage{latexsym}
\usepackage{amssymb}
\usepackage{amsmath,amsfonts,amssymb,theorem,euscript,array,enumerate,amsfonts,mathrsfs}
\usepackage{graphicx}

\newtheorem{Theorem}{Theorem}[part]
\newtheorem{Definition}{Definition}[part]
\newtheorem{Proposition}{Proposition}[part]

\newtheorem{Lemma}{Lemma}[part]
\newtheorem{Corollary}{Corollary}[part]
\newtheorem{Remark}{Remark}[part]

\def\esssup_#1{\underset{#1}{\mathrm{ess\,sup\, }}}
\def\essinf_#1{\underset{#1}{\mathrm{ess\,inf\, }}}

\def \trans{^{\scriptscriptstyle{\intercal}}}

\def \trans{^{\scriptscriptstyle{\intercal }}}

\def \I{\mathbb{I}}
\def \N{\mathbb{N}}
\def \R{\mathbb{R}}

\def \E{\mathbb{E}}
\def \F{\mathbb{F}}

\def \P{\mathbb{P}}

\def \Ac{{\cal A}}
\def \Bc{{\cal B}}
\def \Cc{{\cal C}}

\def \Fc{{\cal F}}

\def \Lc{{\cal L}}

 \def \Nc{{\cal N}}

\def \Tc{{\cal T}}
\def \Uc{{\cal U}}
\def \Vc{{\cal V}}

\def \Vc{{\cal V}}

\def \eps{\varepsilon}

\def \ep{\hbox{ }\hfill$\Box$}

\def\reff#1{{\rm(\ref{#1})}}

\def\beqs{\begin{eqnarray*}}
\def\enqs{\end{eqnarray*}}
\def\beq{\begin{eqnarray}}
\def\enq{\end{eqnarray}}

\allowdisplaybreaks

\begin{document}

\title{Robust feedback switching control:\\ dynamic programming and viscosity solutions\thanks{E. Bayraktar is supported in part by the National Science Foundation under grant  DMS-1613170. H. Pham is  supported in part by FiME (Laboratoire de Finance des March\'es de l'Energie) and
the ``Finance et D\'eveloppement Durable - Approches Quantitatives'' Chair. }}

\author{Erhan BAYRAKTAR\thanks{Department of Mathematics, University of Michigan, 
		\sf erhan@umich.edu} ~~~
		Andrea COSSO\thanks{Laboratoire de Probabilit\'es et Mod\`eles Al\'eatoires, CNRS, UMR 7599, Universit\'e Paris Diderot,
		\sf  cosso@math.univ-paris-diderot.fr}~~~
               Huy{\^e}n PHAM\thanks{Laboratoire de Probabilit\'es et Mod\`eles Al\'eatoires, CNRS, UMR 7599, Universit{\'e} Paris Diderot, and
               CREST-ENSAE,  \sf pham@math.univ-paris-diderot.fr}
             }

\maketitle

\date

\begin{abstract}
We consider a robust switching control problem.  The controller only observes the evolution of the state process, and thus 
uses feedback (closed-loop) switching strategies, a non standard class of switching controls introduced  in this paper.  The adverse player (nature) chooses open-loop controls that represent  the so-called Knightian uncertainty, i.e., misspecifications of the model. The (half) game switcher versus nature is then formulated as  
a two-step (robust)  optimization problem. We develop the stochastic Perron method in this framework, and  prove that it produces a viscosity 
sub and supersolution to a system of  Hamilton-Jacobi-Bellman (HJB)  variational inequalities, which envelope the value function.  Together with a comparison principle, this characterizes the value function of the game as the unique viscosity solution to the HJB equation, and shows as a byproduct the dynamic programming principle for robust feedback switching control problem. 
\end{abstract}

\vspace{3mm}


\vspace{3mm}


\section{Introduction}

Optimal switching is a class of stochastic control problems that has attracted a lot of interest and generated important developments in applied and financial mathematics.  
Swit\-ching control  consists in sequence of interventions that occur at random discrete times due to switching costs, and naturally arises in investment problems  with fixed transaction costs or in real options. The literature on this topic  is quite large and we refer e.g. to 
\cite{TanYon92}, \cite{LyvPha07}, \cite{ElHam09}, \cite{pham09}, \cite{MR2676760}, \cite{MR2653895},  for a treatment by dynamic programming and PDE methods, to \cite{HamZha10}, \cite{HuTan10}, \cite{EliKha14} for the connection with reflected BSDE methods, and to  \cite{DucZer01}, \cite{CarLud10}, \cite{HamJea07} for various applications to  finance and  real options in energy markets.

The standard approach to the study of a switching control problem is to give an evolution for the controlled state process, with assigned drift and diffusion coefficients. These, however, are in practice obtained through estimation procedures and are unlikely to coincide with the \emph{real} coefficients. For this reason, in the present work we study a switching control problem \emph{robust} to a misspecification of the model for the controlled state process.  This is formalized as follows:  Given $s\geq0$, $x\in\R^d$, and a regime $i\in\I_m:=\{1,\ldots,m\}$, let us consider the controlled system of stochastic differential equations, for $t\geq s$:
\begin{equation}
\label{SDE0}
\begin{cases}
X_t \ = \ x + \int_s^t b(X_r,I_r,u_r)dr + \int_s^t\sigma(X_r,I_r,u_r)dW_r, \\
\,\,I_t \ = \ i\,1_{\{s\leq t<\tau_0\}} + \sum_{n\in\N} \iota_n 1_{\{\tau_n\leq t<\tau_{n+1}\}}.
\end{cases}
\end{equation}
The piecewise constant process $I$ denotes the regime value at any time $t$, whose evolution is determined by the controller through the switching control $\alpha=(\tau_n,\iota_n)_{n\in\N}$, while the process $u$, decided by nature,  brings the uncertainty within the model. In the switching control problem with model uncertainty, the objective of the controller is the maximization of the following functional, over a finite time horizon $T$ $<$ 
$\infty$:
\beqs
J(s,x,i;\alpha,u) &:=& \E\bigg[\int_s^T f(X_r^{s,x,i;\alpha,u},I_r^{s,x,i;\alpha,u},u_r)dr + g(X_T^{s,x,i;\alpha,u},I_T^{s,x,i;\alpha,u}) \notag \\
& & \quad - \; \sum_{n\in\N} c(X_{\tau_n}^{s,x,i;\alpha,u},I_{\tau_n^-}^{s,x,i;\alpha,u},I_{\tau_n}^{s,x,i;\alpha,u}) 1_{\{s\leq\tau_n<T\}} \bigg],
\enqs
playing against nature, described by $u$. This leads to the ``robust'' optimization problem
\begin{equation} \label{Two-Step_Opt}
\sup_\alpha\Big(\inf_u J(s,x,i;\alpha,u)\Big).
\end{equation}
What definition and information pattern for the switching control $\alpha$ and for $u$ should we adopt? As a first attempt, if we interpret \eqref{Two-Step_Opt} as a game between the controller and nature, it would be reasonable to formulate it in terms of \emph{nonanticipating strategies against controls}, as in the seminal paper by Elliott \& Kalton \cite{EK72}.  In this case, $\alpha$ is a non-anticipating switching strategy, while $u$ is an open-loop control. Then, the switcher knows the current and past choices made by the opponent (see Section \ref{Strategies} below for more details on this formulation). In the context of robust optimization, the controller does not know in general the choice made by nature. He knows at most the current  state of the system and its past history, that is  the evolution of $X$ and also of $I$ (by keeping track of his previous actions).  For this reason, inspired by \cite{bardi-capuzzo}, \cite{sirbu14b} (see also \cite{krasovski_subbotin} which considers robust controls over feedback strategies in deterministic setting), we take $\alpha$ as a \emph{feedback (also called closed-loop) switching strategy} 
rather than non-anticipating strategy  (namely, we present a \emph{feedback formulation} of a switching control problem, which is quite uncommon in the literature). On the other hand, $u$ can be an \emph{open-loop control} (nature is aware of the all information at disposal). This leads to the formulation of  {\it robust feedback} switching control  problem where  
both players use controls, one in feedback form (the switcher)  and the other in open loop form (the nature), hence different from the Elliott-Kalton formulation  where one player  observes continuously the control (action) 
of the other player. 
  
We develop the stochastic Perron method in this framework of robust feedback switching strategy. This method was initially introduced to analyze linear problems in \cite{BS12}, Dynkin games in \cite{BS14} and \emph{regular} control problems in \cite{BS13}. Later on, it was adapted to analyze exit time problems in \cite{MR3217159}, control problems with state constraints in \cite{2014arXiv1405.4252R}, singular control problems in \cite{2014arXiv1404.7406B}, stochastic differential games in \cite{sirbu14a} and stochastic control with model uncertainty in \cite{sirbu14b}. Stochastic Perron's method is similar to a \emph{verification theorem} and avoids having to go through the dynamic programming principle (DPP) first (which is not known a priori in this context) to show that the value function is a solution to the HJB equation. Actually, the DPP is obtained as a byproduct of the stochastic Perron method  and comparison principle. Unlike the classical verification theorem, the stochastic Perron does not require the a priori smoothness of the value function.  The method is to construct viscosity (semi-) solutions to the HJB equation, which envelope the value function, and relies on the comparison principle of the HJB equation to conclude that the value function is the unique viscosity solution. 
In order to carry out the construction, one needs to define two suitable classes of functions, denoted by $\Vc^-$ and $\Vc^+$, whose elements are known in the literature on stochastic Perron method as stochastic subsolutions ($\Vc^-$) and stochastic supersolutions ($\Vc^+$). The crucial property of $\Vc^-$ and $\Vc^+$ is closedness under minimization/maximization. Moreover, their members stay below/above the value function. 
The technical part of the proof is in showing that the supremum/infimum of the above classes give a viscosity supersolution/subsolution to the HJB equation. 
One of the advantages of the stochastic Perron method is that it allows us to demonstrate that the information available to nature (whether it uses open-loop or feedback strategies) does not affect the value of the game. We do this by constructing the class $\Vc^+$ for an auxiliary problem, whose elements lie by definition above our original value function. Our results here can be thought of as a generalization of the recent work \cite{sirbu14b}, in which the controller uses elementary feedback strategies. In our setting changing the value of control has a switching cost. This changes the nature of the problem as the past action of the controller needs to be stored as a state variable.  
The presence of this additional state variable  brings about several subtle technical issues, which we resolve  in this paper. For example, concatenating the feedback switching strategies need to be done with care (not to incur an additional cost at the time of concatenation), which forces us to make appropriate  changes in defining the class $\Vc^-$.

We should mention that when one can bootstrap the regularity of the viscosity solutions and show that they are classical solutions, one can still use the classical Perron method of Ishii \cite{ishii87}. This program is carried out by \cite{jan_sirbu12} for a stochastic control problem and by \cite{bayraktar_zhang} for a robust stochastic control problem. In general, however, the PDE may not admit a smooth solution and one has to use the generalization of the Perron method, which we called stochastic Perron's method, described above. 
If one attempts to only use the Perron method in \cite{ishii87} to construct viscosity
solutions one faces a major obstacle: without additional knowledge on the properties of value function, it does
not compare with the output of the classical Perron method. In fact this is exactly what happens in \cite{MR2653895}. 
In fact, Section 2 of \cite{MR2653895} shows that the system of variational inequalities has a unique viscosity solution using the classical Perron method. But when they introduce a control problem (not a game) in Section 3, they still go through first proving the DPP, to show that the value function is a viscosity solution and is therefore the unique viscosity solution they constructed in Section 2. 

We should emphasize that although the system of variational inequalities in Section 2 of  \cite{MR2653895}  is quite close to the one in our paper, these authors make the connection in their Section 3 with a control problem for the particular  case when there is one single player using switching and 
regular controls.  Our main result is on one hand  the formulation and solution of  the robust feedback switching control problem, in which the controller only observes the evolution of the state process, and thus uses feedback (closed-loop) switching strategies, a non standard class of switching controls introduced  for the first time in this paper, and on the other hand to prove \emph{directly} that it is the unique viscosity solution to the corresponding system of dynamic programming variational inequalities.

The rest of this paper is organized as follows. In Section \ref{S:ControlProblem}, we provide a rigorous formulation of the robust feedback switching control problem. 
We develop in Section 3 the stochastic Perron method, and characterize the infimum (resp. supremum) of $\Vc^+$ (resp. $\Vc^-$) as the viscosity subsolution (resp. supersolution)
of the HJB equation. In Section 4, by using a comparison principle under a no free loop condition on the switching costs, we conclude that the value function is the unique viscosity solution to the HJB equation, and obtain as a byproduct the dynamic programming principle. We finally compare the two formulations: robust feedback/Elliott-Kalton, in a specific example, which then  gives a counterexample to uniqueness for the HJB equation. In order to keep the paper size reasonable, whenever a result has a standard proof or a similar proof can be found in the literature, we do not report all details, but we focus on the main steps providing a sketch of the proof.

\section{Modeling a robust switching control problem}
\label{S:ControlProblem}

\setcounter{equation}{0}
\setcounter{Theorem}{0} \setcounter{Proposition}{0}
\setcounter{Corollary}{0} \setcounter{Lemma}{0}
\setcounter{Definition}{0} \setcounter{Remark}{0}

\subsection{Feedback switching system under model uncertainty}

In this section, we consider the situation where the switcher knows just  the current and past history of the state.  To model this information pattern, we adopt the notion of {\it feedback strategies}  following the definition introduced 
in the book \cite{bardi-capuzzo} (see Chapter VIII, Section 3.1) or in \cite{sirbu14b}.  It is important to notice that this notion of {\it feedback strategies}  differs from the notion of nonanticipating strategies  \`a la Elliott-Kalton where the 
switcher-player knows the current and past choices of the control made by  his/her  opponent  (here the nature), see also the discussion in  Chapter VIII  of \cite{bardi-capuzzo} and in particular Lemma 3.5 which gives the connection 
between these two notions.

Let $U$ be a compact metric space and $(\Omega,\Fc,\P)$ be a fixed probability space on which a $d$-dimensional Brownian motion $W$ $=$ $(W_t)_{t\geq 0}$ is defined. For any $s\geq0$, we consider the filtration $\F^{W,s}=(\Fc_t^{W,s})_{t\geq s}$, which is the augmented natural filtration generated by the Brownian increments starting at $s$, i.e.,
\beqs
\Fc_t^{W,s} &:=& \sigma(W_r-W_s,\,s\leq r\leq t)\vee\Nc(\P,\Fc), \qquad t\geq s,
\enqs
where $\Nc(\P,\Fc):=\{N\in\Fc\colon\P(N)=0\}$. For each $s\geq0$, we denote by $\F^s=(\Fc_t^s)_{t\geq s}$ another filtration satisfying the usual conditions, which is larger than $\F^{W,s}$ and keeps $(W_t-W_s)_{t\geq s}$ a Brownian motion starting at $s$.

We fix a finite time horizon $0<T<\infty$. For any $s\in[0,T]$, we denote by $y(\cdot)$ or $y$ a generic element of the space $C([s,T];\R^d)\times\mathscr{L}([s,T];\I_m)$, where $\mathscr{L}([s,T];\I_m)$ denotes the set of c\`agl\`ad paths valued in $\I_m$ (notice that the elements of $\mathscr{L}([s,T];\I_m)$ are indeed piecewise constant paths, since $\I_m$ is a discrete set). We also write $y=(y^X,y^I)$ with $y^X\in C([s,T];\R^d)$ and $y^I\in\mathscr{L}([s,T];\I_m)$. We define the filtration $\mathbb B^s=(\Bc_t^s)_{s\leq t\leq T}$, where $\Bc_t^s$ is the $\sigma$-algebra generated by the canonical coordinate maps $C([s,T];\R^d)\times\mathscr{L}([s,T];\I_m)\rightarrow\R^d\times\I_m$, $y(\cdot)\mapsto y(r)$, $r\in[s,t]$, namely
\[
\Bc_t^s \ := \ \sigma(y(\cdot)\mapsto y(r),\,s\leq r\leq t).
\]
A map $\tau\colon C([s,T];\R^d)\times\mathscr{L}([s,T];\I_m)\rightarrow[s,T]$ satisfying $\{\tau\leq t\}\in\Bc_t^s$, $\forall\,t\in[s,T]$, is called a \emph{stopping rule}. $\Tc^s$ denotes the family of all stopping rules starting at $s$. For any $s\in[0,T]$ and $\tau\in\Tc^s$, we define, as usual,
\beqs
\Bc_{\tau^+}^s\!\!\! &:=& \big\{B\in\Bc_T^s\colon\forall\,t\in[s,T],\,B\cap\{y\colon\tau(y)\leq t\}\in\Bc_{t^+}^s\big\}, \\
\Bc_\tau^s &:=& \big\{B\in\Bc_T^s\colon\forall\,t\in[s,T],\,B\cap\{y\colon\tau(y)\leq t\}\in\Bc_t^s\big\},
\enqs
where $\Bc_{t^+}^s:=\cap_{r>t}\Bc_r^s$, $t\in[s,T)$, and $\Bc_{T^+}^s:=\Bc_T^s$. We also denote $y(T^+):=y(T)$, for any $y\in C([s,T];\R^d)\times\mathscr{L}([s,T];\I_m)$.

\begin{Definition}[Feedback switching strategies]
\label{D:switching}
Fix $s\in[0,T]$. We say that the double sequence $\alpha=(\tau_n,\iota_n)_{n\in\N}$ is a feedback switching control starting at $s$ if:
\begin{itemize}
\item $\tau_n\in\Tc^s$, for any $n\in\N$, and
\[
s \ \leq \ \tau_0 \ \leq \ \cdots \ \leq \ \tau_n \ \leq \ \cdots \ \leq \ T.
\]
Moreover, $(\tau_n)_{n\in\N}$ satisfies the following property: $\forall\,(y_n)_{n\in\N}\in C([s,T];\R^d)\times\mathscr{L}([s,T];\I_m)$, with $y_n(t)=y_{n+1}(t)$, $t\in[s,\tau_n(y_n)]$, for every $n\in\N$, then
\[
\tau_n(y_n) \ = \ T, \qquad \text{for $n$ large enough}.
\]
\item $\iota_n\colon C([s,T];\R^d)\times\mathscr{L}([s,T];\I_m)\rightarrow\I_m$ is $\Bc_{\tau_n}^s$-measurable, for any $n\in\N$.
\end{itemize}
$\Ac_s$ denotes the family of all feedback switching controls starting at $s$.
\end{Definition}

\begin{Remark}
{\rm This canonical definition of the feedback switching strategy means that the stopping rules $\tau_n$ are based on the observation of the state, while the actions $\iota_n$ decided at time $\tau_n$  
are based only on the knowledge of the state up to the decision time.   We may alternatively call feedback switching strategy as closed-loop switching control as opposed to the notion of open loop switching controls, where the decision times $\tau_n$ are stopping times with respect to the larger 
filtration $\F^s$, and the actions 
$\iota_n$ are  based on a larger information given by the filtration $\F^s$.  Consider a sequence of paths $(y_n)_{n\in\N}$ as in Definition \ref{D:switching}. Then, the sequence $(\tau_n(y_n))_{n\in\N}$ is nondecreasing. Indeed, from Lemma \ref{L:tau&iota1} below we have $\tau_n(y_n)=\tau_n(y_{n+1})$. Since $\tau_n(y_{n+1})\leq\tau_{n+1}(y_{n+1})$ from the nondecreasing property of the sequence $(\tau_n)_{n\in\N}$, the thesis follows. See also Remark \ref{R:Property_tau_n} below, where the property ``$\tau_n(y_n)=T$, for $n$ large enough'' is analyzed in detail. 
This  structure condition on the sequence $(y_n)$ is required for ensuring well-posedness, i.e. 
 in order to guarantee that the optimal control does not have infinitely many switches and that the SDE \reff{SDE} of X is well defined.  
 This is discussed in detail below, see in particular  Remark 2.3. 
\ep
}
\end{Remark}

\begin{Definition}[Open-loop controls]
Fix $s\in[0,T]$. An open-loop control $u$ starting at $s$, for the nature, is an $\F^s$-progressively measurable process $u\colon[s,T]\times\Omega\rightarrow U$. We denote by $\Uc_s$ the collection of all possible open-loop controls, given the initial deterministic time $s$.
\end{Definition}

For any $(s,x,i)\in[0,T]\times\R^d\times\I_m$, $\alpha=(\tau_n,\iota_n)_{n\in\N}\in\Ac_s$, $u\in\Uc_s$, we can now write equation \eqref{SDE0} on $[0,T]$ as follows:
\begin{equation}
\label{SDE}
\begin{cases}
X_t \ \, \!\!= \ x + \int_s^t b(X_r,I_r,u_r)dr + \int_s^t\sigma(X_r,I_r,u_r)dW_r, &s\leq t\leq T, \\
I_t \ \, = \ i1_{\{s\leq t<\tau_0(X_\cdot,I_{\cdot^-})\}} + \sum_{n\in\N} \iota_n(X_\cdot,I_{\cdot^-}) 1_{\{\tau_n(X_\cdot,I_{\cdot^-})\leq t<\tau_{n+1}(X_\cdot,I_{\cdot^-})\}}, &s\leq t<T, \\
I_T \ \, \!\!= \ I_{T^-},
\end{cases}
\end{equation}
with $I_{s^-}:=I_s$. Notice that the presence of $I_{\cdot^-}$ in place of $I_\cdot$ in the arguments of $\tau_n,\iota_n$ is due to the fact that the choice of $(\tau_n,\iota_n)$ by the controller is based only on the information coming from the previous switching actions $(\tau_i,\iota_i)_{0\leq i\leq n-1}$. Moreover, the last equation $I_T=I_{T^-}$ in \eqref{SDE} means that there is no regime switching at the final time $T$. We impose the following assumptions on the coefficients $b\colon\R^d\times\I_m\times U \rightarrow \R^d$ and $\sigma\colon\R^d\times\I_m\times U \rightarrow \R^{d\times d}$ (in the sequel, we use the notation $\|A\|^2$ $=$ $\text{tr}(AA\trans)$ for the Hilbert-Schmidt norm of any matrix $A$).

\vspace{2mm}

{\bf (H1)}
\begin{itemize}
\item[(i)] $b,\sigma$ are jointly continuous on $\R^d\times\I_m\times U$.
\item[(ii)] $b,\sigma$ are uniformly Lipschitz continuous in $x$, i.e.,
\beqs
|b(x,i,u)-b(x',i,u)| + \|\sigma(x,i,u)-\sigma(x',i,u)\| &\le & L_1 |x-x'|,
\enqs
$\forall\,x,x'\in\R^d$, $i\in\I_m$, $u\in U$, for some positive constant $L_1$.
\end{itemize}

\begin{Remark}
\label{R:LinearGrowth_b_sigma}
{\rm
From Assumption {\bf (H1)} it follows that $b$ and $\sigma$ satisfy a linear growth condition in $x$, i.e.,
\beqs
|b(x,i,u)| + \|\sigma(x,i,u)\| &\leq & M_1(1+|x|),
\enqs
$\forall\,x\in\R^d$, $i\in\I_m$, $u\in U$, for some positive constant $M_1$.
\ep
}
\end{Remark}

\begin{Remark}
\label{R:Property_tau_n}
{\rm
Fix $s\in[0,T]$ and $\alpha=(\tau_n,\iota_n)_{n\in\N}\in\Ac_s$. Let us consider the following properties of the nondecreasing sequence $(\tau_n)_{n\in\N}$:
\begin{itemize}
\item[(i)] \emph{Uniformly finite.} There exists $N\in\N$ such that, $\forall\,y\in C([s,T];\R^d)\times\mathscr{L}([s,T];\I_m)$,
\[
\tau_n(y) \ = \ T, \qquad \text{for $n\geq N$}.
\]
\item[(ii)] \emph{Finite along every adaptive sequence.} For every sequence $(y_n)_{n\in\N}\in C([s,T];\R^d)\times\mathscr{L}([s,T];\I_m)$ satisfying, for every $n\in\N$, $y_n(t)=y_{n+1}(t)$, $\forall\,t\in[s,\tau_n(y_n)]$, we have
\[
\tau_n(y_n) \ = \ T, \qquad \text{for $n$ large enough}.
\]

\item[(iii)] \emph{Finite along every path.} $\forall\,y\in C([s,T];\R^d)\times\mathscr{L}([s,T];\I_m)$,
\[
\tau_n(y) \ = \ T, \qquad \text{for $n$ large enough}.
\]
\end{itemize}
Condition (i) is the strongest, while (iii) is the weakest. In Definition \ref{D:switching} we imposed the intermediate property (ii), since it allows to have a well-posedness result for equation \eqref{SDE}, which is no longer guaranteed if we require only (iii). To see this latter point, we construct a counter-example. Take $s=0$, $T=1$, and $m=2$ so that $\I_2=\{1,2\}$. Consider the sequence $(b_n)_{n\in\N}\subset[0,1]$ given by
\beqs
b_n &=& \sum_{j=0}^n \frac{1}{2^{j+2}}, \qquad \forall\,n\in\N.
\enqs
In particular, we have $b_0=\frac{1}{4}$, $b_1=\frac{1}{4}+\frac{1}{8}$, $b_2=\frac{1}{4}+\frac{1}{8}+\frac{1}{16}$, $\ldots$, and in general
\beqs
b_n &=& \frac{2^{n+1}-1}{2^{n+2}}, \qquad \forall\,n\geq0.
\enqs
Notice that $(b_n)_{n\in\N}$ is a strictly increasing sequence satisfying $b_n\nearrow\frac{1}{2}$, as $n\rightarrow\infty$. Now, for every $y\in C([0,1];\R^d)\times\mathscr{L}([0,1];\I_2)$ we write $y=(y^X,y^I)$ with $y^X\in C([0,1];\R^d)$ and $y^I\in\mathscr{L}([0,1];\I_2)$. Then, we define the sequence $(\tau_n)_{n\in\N}$ as follows:
\beqs
\tau_n(y) &=& b_n 1_{\{y\in B_n\}} + 1_{\{y\in B_n^c\}}, \qquad \forall\,y\in C([0,1];\R^d)\times\mathscr{L}([0,1];\I_2),\,n\in\N,
\enqs
where
\beqs
B_0 &=& \big\{y\in C([0,1];\R^d)\times\mathscr{L}([0,1];\I_2)\colon y^I(t)=y^I(0),\,0< t\leq b_0\big\}, \\
B_n &=& \big\{y\in B_{n-1}\colon y^I(t)=3-y^I(b_{n-1}),\,b_{n-1}< t\leq b_n\big\}, \qquad\qquad\;\;\, \forall\,n\geq1.
\enqs
Observe that, since $y^I(t)\in\I_2$ then $3-y^I(t)\in\I_2$; moreover, when $y^I(t)=1$ then $3-y^I(t)=2$, while if $y^I(t)=2$ then $3-y^I(t)=1$. We also notice that $B_n\in\Bc_{b_n}^0$, therefore $\tau_n\in\Tc^0$. Furthermore, $(\tau_n)_{n\in\N}$ is a nondecreasing sequence which verifies property (iii) above: this is due to the fact that every path $y\in C([0,1];\R^d)\times\mathscr{L}([0,1];\I_2)$ has only a finite number of jumps, since $\I_2$ is a discrete set; in other words, any $y$ belongs to $B_n^c$ when $n$ is large enough (e.g., when $n$ is strictly greater than the number of jumps of $y$). However, $(\tau_n)_{n\in\N}$ does not satisfy property (ii), as we shall prove below. We also define
\beqs
\iota_n(y) &=& 3 - y^I(b_n), \qquad \forall\,y\in C([0,1];\R^d)\times\mathscr{L}([0,1];\I_2),\,n\in\N.
\enqs
In other words, when $y^I(b_n)=1$ then $\iota_n(y)=2$, while when $y^I(b_n)=2$ then $\iota_n(y)=1$. Let $\alpha=(\tau_n,\iota_n)_{n\in\N}$, then $\alpha$ satisfies Definition \ref{D:switching}, but for property (ii) (see below), even if property (iii) is satisfied. Now, we solve equation \eqref{SDE} with $x\in\R^d$, $\alpha=(\tau_n,\iota_n)_{n\in\N}$, $u\in\Uc_{0,0}$, and $i=1\in\I_2$. Define the (deterministic) process $I\colon[0,1]\rightarrow\I_2$ as follows, for any $t\in[0,\frac{1}{2})$,
\beqs
I_t &=&
\begin{cases}
1, \qquad\quad & 0\leq t\leq b_0, \\
2, & b_0< t\leq b_1, \\
1, & b_1< t\leq b_2, \\
2, & b_2< t\leq b_3, \\
\vdots
\end{cases}
\enqs
On the other hand, we do not specify $I$ on $[\frac{1}{2},1]$, we only require that the limit $I_{1^-}:=\lim_{t\uparrow1}I_t$ exists and we suppose that $I_1=I_{1^-}$. Notice that $I_{\frac{1}{2}^-}$ does not exist, therefore $I\notin\mathscr{L}([0,1];\I_2)$. However, the process $I$ solves equation \eqref{SDE} (viceversa, every process satisfying \eqref{SDE} coincides with $I$ on the interval $[0,\frac{1}{2})$; in particular, there does not exist a solution process with paths in $\mathscr{L}([0,1];\I_2)$). Moreover, under Assumption {\bf (H1)} we can also solve equation \eqref{SDE} for $X$. Since we did not specify the behavior of $I$ on the entire interval $[0,1]$, we can not have uniqueness of the solution for \eqref{SDE}. Nevertheless, we notice that the sequence $(\tau_n)_{n\in\N}$ does not satisfy property (ii) above. Indeed, let $y_n(\cdot):=I_{\cdot\wedge b_n}$, $n\in\N$. Then, $y_n\in\mathscr{L}([0,1];\I_2)$, but $\tau_n(y_n)<\frac{1}{2}$, for any $n$. This shows that if we only require property (iii), then the well-posedness of equation \eqref{SDE} is no longer guaranteed.
\ep
}
\end{Remark}

We now study the well-posedness of equation \eqref{SDE}, for which we need the following two lemmata.

\begin{Lemma}
\label{L:tau&iota1}
Let $s\in[0,T]$, $\tau\in\Tc^s$, and $y_1,y_2\in C([s,T];\R^d)\times\mathscr{L}([s,T];\I_m)$. If $y_1(t)=y_2(t)$, $s\leq t\leq\tau(y^1)$, then:
\begin{itemize}
\item[\textup{(i)}] $\tau(y_1)=\tau(y_2)$.
\item[\textup{(ii)}] $\iota(y_1)=\iota(y_2)$, for any $\Bc_\tau^s$-measurable map $\iota\colon C([s,T];\R^d)\times\mathscr{L}([s,T];\I_m)\rightarrow\I_m$.
\end{itemize}
\end{Lemma}
\textbf{Proof}
Let $t^*:=\tau(y_1)$. We begin noting that if $B\in\Bc_{t^*}^s$ and $y_1\in B$, then $y_2\in B$, as well. Since $\tau$ is a stopping rule, the event $B:=\{y\colon\tau(y)=t^*\}$ belongs to $\Bc_{t^*}^s$. As $y_1\in B$, we then see that $y_2\in B$, i.e., $\tau(y_2)=\tau(y_1)$, which gives (i). Notice that assertion (i) can be also deduced by (100.1) at page 149, Chapter IV, in \cite{dellacherie_meyerI}.

Concerning (ii), let $\iota\colon C([s,T];\R^d)\times\mathscr{L}([s,T];\I_m)\rightarrow\I_m$ be $\Bc_\tau^s$-measurable. By definition of $\iota$, the event $\tilde B:=\{y\colon\iota(y)=\iota(y_1)\}$ belongs to $\Bc_\tau^s$. Therefore, $B:=\tilde B\cap\{\tau(y)\leq t^*\}\in\Bc_{t^*}^s$. Since $y_1\in B$, from the observation at the beginning of the proof it follows that $y_2\in B$, which implies $y_2\in\tilde B$, i.e., $\iota(y_2)=\iota(y_1)$.
\ep

\begin{Lemma}
\label{L:tau&iota2}
Let $s\in[0,T]$, $\tau\in\Tc^s$, and $Y=(Y_t)_{s\leq t\leq T}$ be an $\F^s$-adapted process valued in $\R^d\times\I_m$. Suppose that every path of $Y$ belongs to $C([s,T];\R^d)\times\mathscr{L}([s,T];\I_m)$. Then, $\tau_Y\colon\Omega\rightarrow[s,T]$ defined as $\tau_Y(\omega):=\tau(Y_\cdot(\omega))$, $\omega\in\Omega$, is an $\F^s$-stopping time. Moreover, if $\iota\colon C([s,T];\R^d)\times\mathscr{L}([s,T];\I_m)\rightarrow\I_m$ is $\Bc_\tau^s$-measurable then $i_Y(\omega):=\iota(Y_\cdot(\omega))$, $\omega\in\Omega$, is $\Fc_{\tau_Y}^s$-measurable.
\end{Lemma}
\textbf{Proof.}
For any $t\in[s,T]$, we notice that the map $Y_\cdot$ is measurable from $(\Omega,\Fc_t^s)$ into $(C([s,T];\R^d)\times\mathscr{L}([s,T];\I_m),\Bc_t^s)$. Then, $\{\omega\colon\tau_Y(\omega)\leq t\}=\{\omega\colon\tau(Y_\cdot(\omega))\leq t\}=\{\omega\colon Y_\cdot(\omega)\in\tau^{-1}([s,t])\}$. Since $\tau^{-1}([s,t])\in\Bc_t^s$, we have $\{\omega\colon Y_\cdot(\omega)\in\tau^{-1}([s,t])\}\in\Fc_t^s$, which implies that $\tau_Y$ is an $\F^s$-stopping time.

Let now $\iota\colon C([s,T];\R^d)\times\mathscr{L}([s,T];\I_m)\rightarrow\I_m$ be $\Bc_\tau^s$-measurable. We have to prove that $\{\omega\colon\iota_Y(\omega)=\underline i\}\in\Fc_{\tau_Y}^s$, for any $\underline i\in\I_m$, i.e., $\{\omega\colon\iota_Y(\omega)=\underline i\}\cap\{\omega\colon\tau_Y(\omega)\leq t\}\in\Fc_t^s$, for any $\underline i\in\I_m$ and $t\in[s,T]$. Then, fix $\underline i\in\I_m$ and $t\in[s,T]$. We have
\begin{align*}
\big\{\omega\colon\iota_Y(\omega)=\underline i\big\}\cap\big\{\omega\colon\tau_Y(\omega)\leq t\big\} \ &= \ \big\{\omega\colon Y_\cdot(\omega)\in\iota^{-1}(\underline i)\big\}\cap\big\{\omega\colon Y_\cdot(\omega)\in\tau^{-1}([s,t])\big\} \\
&= \ \big\{\omega\colon Y_\cdot(\omega)\in\{y\colon\iota(y)=\underline i\}\cap\{y\colon\tau(y)\leq t\}\big\}.
\end{align*}
Since $\iota$ is $\Bc_\tau^s$-measurable, then $\{y\colon\iota(y)=\underline i\}\cap\{y\colon\tau(y)\leq t\}\in\Bc_t^s$. Therefore, from the observation at the beginning of the proof, we get the thesis.
\ep

\begin{Proposition}
\label{P:X}
Let Assumption {\bf (H1)} hold. For any $(s,x,i)\in[0,T]\times\R^d\times\I_m$, $\alpha\in\Ac_s$, $u\in\Uc_s$, there exists a unique (up to indistinguishability) $\F^s$-adapted process $(X^{s,x,i;\alpha,u},I^{s,x,i;\alpha,u})=(X_t^{s,x,i;\alpha,u},I_t^{s,x,i;\alpha,u})_{s\leq t\leq T}$ to equation \eqref{SDE}, such that every path of $(X_\cdot^{s,x,i;\alpha,u},I_{\cdot^-}^{s,x,i;\alpha,u})$ belongs to $C([s,T];\R^d)\times\mathscr{L}([s,T];\I_m)$. Moreover, for any $q\geq1$ there exists a positive constant $C_{q,T}$, depending only on $q,T,M_1$ (independent of $s,x,i,\alpha,u$), such that
\beq
\label{EstimateX}
\E\Big[\sup_{s\leq t\leq T}|X_t^{s,x,i;\alpha,u}|^q\Big] &\leq & C_{q,T}(1+|x|^q).
\enq
\end{Proposition}
\begin{Remark}
{\rm
In Proposition \ref{P:X} we require that every path of $(X_\cdot^{s,x,i;\alpha,u},I_{\cdot^-}^{s,x,i;\alpha,u})$ belongs to $C([s,T];\R^d)\times\mathscr{L}([s,T];\I_m)$ in order to guarantee that the maps $\tau_n(X_\cdot^{s,x,i;\alpha,u}(\omega),I_{\cdot^-}^{s,x,i;\alpha,u}(\omega))$ and $\iota_n(X_\cdot^{s,x,i;\alpha,u}(\omega),I_{\cdot^-}^{s,x,i;\alpha,u}(\omega))$ are well-defined for every $\omega\in\Omega$, $n\in\N$.
\ep
}
\end{Remark}
\textbf{Proof.}
Fix $(s,x,i)\in[0,T]\times\R^d\times\I_m$, $\alpha=(\tau_n,\iota_n)_{n\in\N}\in\Ac_s$, $u\in\Uc_s$.

\noindent\textbf{Step I}. \emph{Existence.} We begin noting that, since the control $\alpha$ is of feedback type, we have to construct the solution $(X^{s,x,i;\alpha,u},I^{s,x,i;\alpha,u})$ and $\alpha$ simultaneously. To do it we proceed as follows: for any $N\in\N$, we solve equation \eqref{SDE} controlled by $u$ and the first $N$ switching actions $(\tau_n,\iota_n)_{0\leq n\leq N-1}$. This is done by induction on $N$. Then, noting that $(X^N,I^N)=(X^{N-1},I^{N-1})$ on the stochastic interval $[s,\tau_{N-1})$, by pasting together the various solutions we are able to construct a solution $(X^{s,x,i;\alpha,u},I^{s,x,i;\alpha,u})$ to the original equation \eqref{SDE} with the entire switching control $\alpha$. We now report the rigorous arguments.

For any $N\in\N$, let $\alpha^N=(\tau_n^N,\iota_n^N)_{n\in\N}\in\Ac_s$ be given by
\beqs
(\tau_n^N,\iota_n^N) &:=&
\begin{cases}
(\tau_n,\iota_n), \qquad & 0 \leq n \leq N-1, \\
(T,\iota_n), & n \geq N.
\end{cases}
\enqs
Let $N=0$ and consider equation \eqref{SDE} controlled by $\alpha^0$ and $u$. Notice that $I$ is uncontrolled, in particular $I_t=i$, $s\leq t\leq T$. Then, it is well-known that under Assumption {\bf (H1)} there exists a unique (up to indistinguishability) $\F^s$-adapted solution $(X_t^0,I_t^0)_{s\leq t\leq T}$ to this equation, with  $I_t^0=i$ for any $t\in[s,T]$, such that every (not only $\P$-a.e., simply choosing an opportune indistinguishable version) path of $(X_\cdot^0,I_{\cdot^-}^0)$ belongs to $C([s,T];\R^d)\times\mathscr{L}([s,T];\I_m)$.

Now, let us prove the inductive step. Let $N\in\N\backslash\{0\}$ and suppose that there exists an $\F^s$-adapted solution $(X^{N-1},I^{N-1})$ to equation \eqref{SDE} controlled by $\alpha^{N-1}$ and $u$, such that every path of $(X_\cdot^{N-1},I_{\cdot^-}^{N-1})$ belongs to $C([s,T];\R^d)\times\mathscr{L}([s,T];\I_m)$. Our aim is to solve equation \eqref{SDE} controlled by $\alpha^N$ and $u$. To this end, we define the process $I^N=(I_t^N)_{s\leq t\leq T}$ as follows:
\beqs
\begin{cases}
I_t^N \ = \ I_t^{N-1} 1_{\{s\leq t<\tau_{N-1}(X_\cdot^{N-1},I_{\cdot^-}^{N-1})\}} + \iota_{N-1}(X_\cdot^{N-1},I_{\cdot^-}^{N-1}) 1_{\{\tau_{N-1}(X_\cdot^{N-1},I_{\cdot^-}^{N-1})\leq t<T\}}, \\
I_T^N \ = \ I_{T^-}^N.
\end{cases}
\enqs
From Lemma \ref{L:tau&iota2} we see that $I^N$ is an $\F^s$-adapted process, with every path in $\mathscr{L}([s,T];\I_m)$. Then, under Assumption {\bf (H1)} there exists a unique (up to indistinguishability) $\F^s$-adapted solution $(X_t^N,I_t^N)_{s\leq t\leq T}$ to equation \eqref{SDE}, such that every path of $(X_\cdot^N,I_{\cdot^-}^N)$ belongs to $C([s,T];\R^d)\times\mathscr{L}([s,T];\I_m)$. Since $(X^N,I^N)$ and $(X^{N-1},I^{N-1})$ solve the same equation on $[s,\tau_{N-1}(X_\cdot^{N-1},I_{\cdot^-}^{N-1}))$, then $(X_t^N,I_t^N)=(X_t^{N-1},I_t^{N-1})$, $t\in[s,\tau_{N-1}(X_\cdot^{N-1},I_{\cdot^-}^{N-1}))$. In particular, $(X_t^N,I_{t^-}^N)=(X_t^{N-1},I_{t^-}^{N-1})$, for any $t\in[s,\tau_{N-1}(X_\cdot^{N-1},I_{\cdot^-}^{N-1})]$. From Lemma \ref{L:tau&iota1}, it follows that
\beqs
\big(\tau_n(X_\cdot^{N-1},I_{\cdot^-}^{N-1}),\iota_n(X_\cdot^{N-1},I_{\cdot^-}^{N-1})\big) &=& \big(\tau_n(X_\cdot^N,I_{\cdot^-}^N),\iota_n(X_\cdot^N,I_{\cdot^-}^N)\big), \qquad 0\leq n\leq N-1.
\enqs
As a consequence, $(X^N,I^N)$ solves equation \eqref{SDE} controlled by $\alpha^N$ and $u$.

Finally, let us define (with the convention $\tau_{-1}:=s$)
\beq
X_t^{s,x,i;\alpha,u} &:=& \sum_{n\in\N} X_t^N 1_{\{\tau_{N-1}(X_\cdot^{N-1},I_{\cdot^-}^{N-1})\leq t<\tau_N(X_\cdot^N,I_{\cdot^-}^N)\}}, \label{X} \\
I_t^{s,x,i;\alpha,u} &:=& \sum_{n\in\N} I_t^N 1_{\{\tau_{N-1}(X_\cdot^{N-1},I_{\cdot^-}^{N-1})\leq t<\tau_N(X_\cdot^N,I_{\cdot^-}^N)\}}, \label{I}
\enq
for any $s\leq t< T$ and $(X_T^{s,x,i;\alpha,u},I_T^{s,x,i;\alpha,u}):=(X_{T^-}^{s,x,i;\alpha,u},I_{T^-}^{s,x,i;\alpha,u})$. For simplicity of notation, we denote $(X,I):=(X^{s,x,i;\alpha,u},I^{s,x,i;\alpha,u})$. Recalling that $\tau_{N-1}(X_\cdot^{N-1},I_{\cdot^-}^{N-1})=\tau_{N-1}(X_\cdot^N,I_{\cdot^-}^N)\leq\tau_N(X_\cdot^N,I_{\cdot^-}^N)$, we see that the sequence $(\tau_N(X_\cdot^N,I_{\cdot^-}^N))_{N\geq-1}$ is nondecreasing, so that, for any $t\in[s,T]$, there is at most one term different from zero in the series appearing in \eqref{X} and \eqref{I}. Moreover, from Definition \ref{D:switching}, and, more precisely, from property (ii) of Remark \ref{R:Property_tau_n}, we have that, for every $\omega\in\Omega$, $\tau_N(X_\cdot^N(\omega),I_{\cdot^-}^N(\omega))=T$, for $N$ large enough. In particular, $X$ and $I$ are well-defined over the entire interval $[s,T]$ and they are $\F^s$-adapted. Furthermore, we notice that $(X_t,I_t)=(X_t^N,I_t^N)$, $t\in[s,\tau_N(X_\cdot^N,I_{\cdot^-}^N))$. Then, using again property (ii) of Remark \ref{R:Property_tau_n}, it follows that every path of $(X_\cdot,I_{\cdot^-})$ belongs to $C([s,T];\R^d)\times\mathscr{L}([s,T];\I_m)$. In addition, since $(X_t,I_{t^-})=(X_t^N,I_{t^-}^N)$, $t\in[s,\tau_N(X_\cdot^N,I_{\cdot^-}^N)]$, from Lemma \ref{L:tau&iota1} we have
\beqs
\big(\tau_N(X_\cdot^N,I_{\cdot^-}^N),\iota_N(X_\cdot^N,I_{\cdot^-}^N)\big) &=& \big(\tau_N(X_\cdot,I_{\cdot^-}),\iota_N(X_\cdot,I_{\cdot^-})\big), \qquad \forall\,N\in\N.
\enqs
In particular, $(X_t,I_t)=(X_t^N,I_t^N)$, $t\in[s,\tau_N(X_\cdot,I_{\cdot^-}))$. This implies that $(X,I)$ solves equation \eqref{SDE} on $[s,\tau_N(X_\cdot,I_{\cdot^-}))$, for any $N\in\N$. Recalling property (ii) of Remark \ref{R:Property_tau_n}, we see that $(X,I)$ solves equation \eqref{SDE} on $[s,T)$. Since, by definition, $(X_T,I_T)=(X_{T^-},I_{T^-})$, it follows that $(X,I)$ solves equation \eqref{SDE} on $[s,T]$.

\noindent\textbf{Step II}. \emph{Uniqueness.} Let $(X^1,I^1)$ and $(X^2,I^2)$ be two solutions of \eqref{SDE}. Set $\underline\tau_0:=\tau_0(X_\cdot^1,I_{\cdot^-}^1)\wedge\tau_0(X_\cdot^2,I_{\cdot^-}^2)$. Notice that $(X^1,I^1)$ and $(X^2,I^2)$ solve the same equation on $[0,\underline\tau_0)$. Therefore $(X^1,I^1)$ and $(X^2,I^2)$ are equal (up to indistinguishability) on $[0,\underline\tau_0)$. Consider $\omega\in\Omega$ such that $\underline\tau_0(\omega)=\tau_0(X_\cdot^1(\omega),I_{\cdot^-}^1(\omega))$. Since $(X_t^1(\omega),I_{t^-}^1(\omega))=(X_t^2(\omega),I_{t^-}^2(\omega))$, $t\in[s,\underline\tau_0(\omega)]=[s,\tau_0(X_\cdot^1(\omega),I_{\cdot^-}^1(\omega))]$, from Lemma \ref{L:tau&iota1} it follows that $\tau_0(X_\cdot^1(\omega),I_{\cdot^-}^1(\omega))=\tau_0(X_\cdot^2(\omega),I_{\cdot^-}^2(\omega))$. When $\underline\tau_0(\omega)=\tau_0(X_\cdot^2(\omega),I_{\cdot^-}^2(\omega))$, a similar argument shows that we still have $\tau_0(X_\cdot^1(\omega),I_{\cdot^-}^1(\omega))=\tau_0(X_\cdot^2(\omega),I_{\cdot^-}^2(\omega))$. From the arbitrariness of $\omega$, we conclude that $\underline\tau_0=\tau_0(X_\cdot^1,I_{\cdot^-}^1)=\tau_0(X_\cdot^2,I_{\cdot^-}^2)$. Using again Lemma \ref{L:tau&iota1}, we also deduce $\iota_0(X_\cdot^1,I_{\cdot^-}^1)=\iota_0(X_\cdot^2,I_{\cdot^-}^2)$. By induction on $n$, we can prove that
\beqs
\big(\tau_n(X_\cdot^1,I_{\cdot^-}^1),\iota_n(X_\cdot^1,I_{\cdot^-}^1)\big) &=& \big(\tau_n(X_\cdot^2,I_{\cdot^-}^2),\iota_n(X_\cdot^2,I_{\cdot^-}^2)\big), \qquad \forall\,n\in\N, \\
(X_t^1,I_t^1) &=& (X_t^2,I_t^2), \hspace{3.6cm} \forall\,t\in[s,\tau_n(X_\cdot^1,I_{\cdot^-}^1)),\,n\in\N.
\enqs
From Definition \ref{D:switching}, and, more precisely, from property (ii) of Remark \ref{R:Property_tau_n}, we have that, for any $\omega\in\Omega$, $\tau_n(X_\cdot^1(\omega),I_{\cdot^-}^1(\omega))=T$ for $n$ large enough. As a consequence, $(X^1,I^1)$ and $(X^2,I^2)$ are equal (up to indistinguishability) on $[s,T)$. Since $(X_T^1,I_T^1)=(X_{T^-}^1,I_{T^-}^1)$ and $(X_T^2,I_T^2)=(X_{T^-}^2,I_{T^-}^2)$, we conclude that $(X^1,I^1)$ and $(X^2,I^2)$ are equal (up to indistinguishability) on $[s,T]$.

\noindent\textbf{Step III}. \emph{Estimate \eqref{EstimateX}.} Under {\bf (H1)}, estimate \eqref{EstimateX} is well-known, see, e.g., Theorem 1.3.15 in \cite{pham09}.
\ep

 \begin{Remark}
{\rm 
Notice  that $\F^s$  is the filtration generated by the noise and $\mathbb B^s$ is the filtration generated by the state variable $X$. 
Since we have strong existence the latter is a subset of the former but not vice versa since the volatility is allowed to degenerate. $\alpha$ is the control of the switcher (the maximizer of our problem) and it is of feedback type. That is the switcher is only allowed to make a decision by observing the state variable. He is not allowed to observe the noise or the actions of the nature, which uses open loop control, i.e., its control is adapted to $\F^s$. 
\ep
}
\end{Remark}

\subsection{The Value function}

The value function associated to the robust switching control problem is defined as follows:
\beq
\label{V}
V(s,x,i) &:=& \sup_{\alpha\in\Ac_s}\inf_{u\in\Uc_s} J(s,x,i;\alpha,u), \qquad \forall\,(s,x,i)\in[0,T]\times\R^d\times\I_m,
\enq
with
\beq
J(s,x,i;\alpha,u) &:=& \E\bigg[\int_s^T f(X_r^{s,x,i;\alpha,u},I_r^{s,x,i;\alpha,u},u_r)dr + g(X_T^{s,x,i;\alpha,u},I_T^{s,x,i;\alpha,u}) \notag \\
& & \quad - \; \sum_{n\in\N} c(X_{\tau_n}^{s,x,i;\alpha,u},I_{\tau_n^-}^{s,x,i;\alpha,u},I_{\tau_n}^{s,x,i;\alpha,u}) 1_{\{s\leq\tau_n<T\}} \bigg], \label{J}
\enq
where $\tau^n$ stands for $\tau^n(X_\cdot^{s,x,i;\alpha,u},I_{\cdot^-}^{s,x,i;\alpha,u})$.

\begin{Remark}
{\rm 
This definition of game value function with the outside player (switcher) using feedback strategies (i.e. closed loop controls) and the inside player (nature) using open-loop controls is 
the same than the one used in Definition 3.6,  Chapter VIII in \cite{bardi-capuzzo}, and called there $B$-feedback value.  It is also pointed out that the $B$-feedback value  is smaller than the upper value of a game 
where the outside player uses nonanticipating strategies \`a la Elliott-Kalton, see also our Section 4.2.   
\ep
}
\end{Remark}

We impose the following conditions on the functions $g\colon\R^d\times\I_m\rightarrow\R$, $f\colon\R^d\times\I_m\times U\rightarrow\R$, and $c\colon\R^d\times\I_m\times\I_m\rightarrow\R$.

\vspace{2mm}

{\bf (H2)}
\begin{itemize}
\item[(i)] $g,f,c$ are jointly continuous on their domains.
\item[(ii)] $c$ is nonnegative.
\item[(iii)] $g,f,c$ satisfy a polynomial growth condition in $x$, i.e.,
\beqs
|g(x,i)| + |f(x,i,u)| + |c(x,i,j)| &\leq & M_2(1+|x|^p),
\enqs
$\forall\,x\in\R^d$, $i,j\in\I_m$, $u\in U$, for some positive constants $M_2$ and $p\geq1$.
\item[(iv)] $g$ satisfies
\beqs
g(x,i) &\geq & \max_{j\neq i}\big[g(x,j) - c(x,i,j)\big],
\enqs
for any $x\in\R^d$ and $i\in\I_m$.
\end{itemize}

\begin{Remark}
\label{R:V_polynomial_growth}
{\rm
Notice that 
$V$ 
satisfies the polynomial growth condition:
\beq
\label{V_linear_growth}
|V(s,x,i)| &\leq & C(1+|x|^p), \qquad \forall\,(s,x,i)\in[0,T]\times\R^d\times\I_m,
\enq
for some positive constant $C$, depending only on $T,M_1,M_2$, and with the same $p$ as in Assumption {\bf (H2)}(iii). Indeed, since $c$ is nonnegative, we find
\beq
\label{UpperBoundV_i}
V(s,x,i) &\leq & \sup_{\alpha\in\Ac_s}\inf_{u\in\Uc_s} \E\bigg[\int_s^T f(X_r^{s,x,i;\alpha,u},I_r^{s,x,i;\alpha,u},u_r)dr + g(X_T^{s,x,i;\alpha,u})\bigg].
\enq
On the other hand, let $\alpha^*=(\tau_n^*,\iota_n^*)_{n\in\N}\in\Ac_s$ be given by $(\tau_n^*,\iota_n^*)=(T,\underline i)$, $\forall\,n\in\N$, for some fixed $\underline i\in\I_m$. Then
\beq
\label{LowerBoundV_i}
V(s,x,i) &\geq & \inf_{u\in\Uc_s} J(s,x,i;\alpha^*,u) \notag \\
&=& \inf_{u\in\Uc_s} \E\bigg[\int_s^T f(X_r^{s,x,i;\alpha^*,u},I_r^{s,x,i;\alpha^*,u},u_r)dr + g(X_T^{s,x,i;\alpha^*,u})\bigg].
\enq
From \eqref{UpperBoundV_i} and \eqref{LowerBoundV_i}, we obtain
\beqs
|V(s,x,i)| &\leq & \sup_{\alpha\in\Ac_s}\sup_{u\in\Uc_s} \E\bigg[\int_s^T |f(X_r^{s,x,i;\alpha,u},I_r^{s,x,i;\alpha,u},u_r)|dr + |g(X_T^{s,x,i;\alpha,u})|\bigg].
\enqs
Now, from estimate \eqref{EstimateX} and the polynomial growth condition of $f$ and $g$ in {\bf (H2)}(iii), we see that estimate \eqref{V_linear_growth} holds. As a consequence, in \eqref{V} we could take the supremum only over $\alpha=(\tau_n,\iota_n)_{n\in\N}\in\Ac_s$ satisfying ($\tau^n$ stands for $\tau^n(X_\cdot^{s,x,i;\alpha,u},I_{\cdot^-}^{s,x,i;\alpha,u})$)
\beqs
\inf_{u\in\Uc_s}\E\bigg[-\sum_{n\in\N} c(X_{\tau_n}^{s,x,i;\alpha,u},I_{\tau_n^-}^{s,x,i;\alpha,u},I_{\tau_n}^{s,x,i;\alpha,u}) 1_{\{s\leq\tau_n<T\}} \bigg] &>& -\infty.
\enqs
\ep
}
\end{Remark}

Our aim is to prove that $V$ is the unique viscosity solution to the dynamic programming equation associated to the robust switching control problem, which turns out to be a system of variational inequalities of Hamilton-Jacobi-Bellman type of the following form:
\beq
\label{HJB}
\begin{cases}
\min\Big\{-\dfrac{\partial V}{\partial t}(s,x,i) - \inf_{u\in U}\big[\Lc^{i,u}V(s,x,i) + f(x,i,u)\big], \\
\hspace{4mm} V(s,x,i) - \max_{j\neq i}\big[V(s,x,j) - c(x,i,j)\big]\Big\} \ = \ 0, \quad (s,x,i)\in[0,T)\times\R^d\times\I_m, \\
V(T,x,i) \ = \ g(x,i), \quad (x,i)\in\R^d\times\I_m,
\end{cases}
\enq
where
\beqs
\Lc^{i,u}V(s,x,i) &=& b(x,i,u).D_x V(s,x,i) + \frac{1}{2}\text{tr}\big[\sigma\sigma\trans(x,i,u)D_x^2 V(s,x,i)\big].
\enqs

We need the definition of (discontinuous) viscosity solution to equation \eqref{HJB}, that we now provide. To this end, given a locally bounded function $v\colon[0,T)\times \R^d\times\I_m\rightarrow\R$, we define its lower semicontinuous (lsc for short) envelope $v_*\colon[0,T]\times \R^d\times\I_m\rightarrow\R$, and upper semicontinuous (usc for short) envelope $v^*\colon[0,T]\times \R^d\times\I_m\rightarrow\R$, by
\beqs
v_*(s,x,i) \;\, = \!\liminf_{\substack{(s',x')\rightarrow (s,x) \\ (s',x')\in[0,T)\times\R^d}} \!\!\!v(s',x',i) \quad\; \text{ and } \quad\; v^*(s,x,i) \;\, = \!\limsup_{\substack{(s',x')\rightarrow (s,x) \\ (s',x')\in[0,T)\times\R^d}} \!\!\!v(s',x',i),
\enqs
for all $(s,x,i)\in[0,T]\times\R^d\times\I_m$.

\begin{Definition}[Viscosity solution to \eqref{HJB}]
\label{D:Viscosity}
\quad \\
\textup{(i)} A lsc $($resp. usc$)$ function $v$ on $[0,T]\times\R^d\times\I_m$ is called a viscosity supersolution $($resp. subsolution$)$ to \eqref{HJB} if
\beqs
v(T,x,i) \ \geq \ (resp. \; \leq) \ g(x,i)
\enqs
for any $(x,i)\in\R^d\times\I_m$, and
\beqs
\min\Big\{-\dfrac{\partial \varphi}{\partial t}(s,x) - \inf_{u\in U}\big[\Lc^{i,u}\varphi(s,x) + f(x,i,u)\big], & & \\
v(s,x,i) - \max_{j\neq i}\big[v(s,x,j) - c(x,i,j)\big]\Big\} &\geq & (resp. \;\, \leq) \;\,  0
\enqs
for any $(s,x,i)\in[0,T)\times\R^d\times\I_m$ and any $\varphi\in C^{1,2}([0,T]\times\R^d)$ such that
\beqs
v(s,x,i) - \varphi(s,x) &=& \min_{(s',x')\in[0,T]\times\R^d} \big[v(s',x',i) - \varphi(s',x')\big] \\
 \Big\{resp. \quad v(s,x,i) - \varphi(s,x) &=& \max_{(s',x')\in[0,T]\times\R^d} \big[v(s',x',i) - \varphi(s',x')\big]\Big\}.
\enqs
\textup{(ii)} A locally bounded function $v$ on $[0,T)\times\R^d\times\I_m$ is called a viscosity solution to \eqref{HJB} if $v_*$ is a viscosity supersolution and $v^*$ is a viscosity subsolution to \eqref{HJB}.
\end{Definition}

\section{Stochastic Perron's method}

\setcounter{equation}{0}
\setcounter{Theorem}{0} \setcounter{Proposition}{0}
\setcounter{Corollary}{0} \setcounter{Lemma}{0}
\setcounter{Definition}{0} \setcounter{Remark}{0}

Our aim is to prove that $V$ is a viscosity solution to the dynamic programming equation \eqref{HJB} and satisfies the dynamic programming principle. To derive these results, we exploit \emph{stochastic Perron's method}, which allows to obtain the viscosity properties of $V$ without relying on the dynamic programming principle, but by means of the comparison theorem for viscosity solutions to \eqref{HJB} (the dynamic programming principle will be obtained as a by-product of this procedure).

\subsection{An Auxiliary robust switching problem}

We begin with the formulation of an auxiliary robust switching control problem where nature adopts closed-loop controls (also called feedback strategies)  in place of open-loop controls. Using the comparison principle for equation \eqref{HJB}, we shall see that the corresponding value function, denoted by $\overline V$, coincides with $V$. In other words, the information available to nature does not affect the value of the game. This is not the only motivation for the introduction of this auxiliary robust control problem. Indeed, in the implementation of the stochastic Perron method we encountered the following difficulty: given two different controls $u_1$ and $u_2$, for nature, we have to concatenate them at some stopping rule $\tau=\tau(X_\cdot,I_{\cdot^-})$. If $u^1$ and $u^2$ are open-loop controls, the control $u^1\otimes_\tau u^2$ resulting from the concatenation of $u^1$ and $u^2$ at the stopping rule $\tau$, given by
\beqs
(u^1\otimes_\tau u^2)(t,\omega,y) &=& u^1(t,\omega)1_{\{s\leq t\leq\tau(y)\}} + u^2(t,\omega)1_{\{\tau(y)<t\leq T\}},
\enqs
is no more of open-loop type, since it also depends on $y$. On the other hand, if $u^1$ and $u^2$ are closed-loop controls, then $u^1\otimes_\tau u^2$ is still a closed-loop control. For this technical reason, to study the original control problem with corresponding value function $V$, we also need to consider another robust switching control problem, in which nature adopts closed-loop controls. In particular, inspired by \cite{sirbu14a} and \cite{sirbu14b}, it turns out that it is more convenient, and it is enough, to consider only piecewise constant closed-loop controls, i.e., the \emph{elementary feedback strategies} that we now define.

\begin{Definition}[Elementary feedback strategies]
\label{D:ElemStrategy}
Fix $s\in[0,T]$. We say that $u$ is an elementary feedback strategy starting at $s$ if:
\begin{itemize}
\item $\tau_k\in\Tc^s$, for any $k=1,\ldots,n$, and
\[
s \ =: \ \tau_0 \ \leq \ \cdots \ \leq \ \tau_k \ \leq \ \cdots \ \leq \ \tau_n \ = \ T.
\]
\item $\xi_k\colon C([s,T];\R^d)\times\mathscr{L}([s,T];\I_m)\rightarrow U$ is $\Bc_{\tau_{k-1}^+}^s$-measurable, for any $k=1,\ldots,n$.
\end{itemize}
The control  $u\colon[s,T]\times C([s,T];\R^d)\times\mathscr{L}([s,T];\I_m)\rightarrow U$ is given by
\[
u(t,y) \ := \ \xi_1(y)1_{\{t=s\}} + \sum_{k=1}^n \xi_k(y)1_{\{\tau_{k-1}(y)<t\leq\tau_k(y)\}}.
\]
$\Uc_s^E$ denotes the family of all elementary feedback strategies (also called elementary closed loop controls) starting at $s$.
\end{Definition}

\begin{Remark}
\label{R:Strategy_tau^+}
{\rm
We notice that Definition \ref{D:ElemStrategy} is inspired by Definition 2.2 in \cite{sirbu14a} (see also Definition 2.1 in \cite{sirbu14b}), the only difference being that $\xi_k$ is $\Bc_{\tau_{k-1}^+}^s$-measurable instead of $\Bc_{\tau_{k-1}}^s$-measurable. This implies that the map $\xi_k=\xi_k(y)$ depends on $y$ through the values $\{y(t),\,s\leq t\leq \tau_{k-1}(y)\}\cup\{y(\tau_{k-1}(y)^+)\}$, so that $\xi_k$ can also depend on $y(\tau_{k-1}(y)^+)$. Recalling that in our setting $y$ denotes a generic path of $(X_t,I_{t^-})_{s\leq t\leq T}$, this means that $\xi_k$ depends on $(X_t,I_t)_{s\leq t\leq\tau_{k-1}(X_\cdot,I_{\cdot^-})}$ rather than on $(X_t,I_{t^-})_{s\leq t\leq\tau_k(X_\cdot,I_{\cdot^-})}$. Therefore, nature reacts to the switcher using all the information at disposal at time $\tau_{k-1}=\tau_{k-1}(X_\cdot,I_{\cdot^-})$, including $I_{\tau_{k-1}}$ (in particular, if $\tau_{k-1}$ coincides with a switching action, nature is aware of the action that the switcher has just performed). We point out that elementary feedback strategies are different from strategies in the sense of Elliott-Kalton where  strategies are used by the outside player (i.e. the switcher here) and not by 
the inside player (the nature here).  Actually, the set of elementary feedback strategies (closed-loop controls) is obviously a subset of open loop controls since they correspond to controls which are piecewise constant on one hand, and with actions decided based only on the knowledge of the state, hence with less information than the one generated by $\F^s$.  In other words, we have $\Uc_s^E$ $\subset$ $\Uc_s$:  for any feedback control $u \in  \mathcal{U}^E_s$ we can construct an open loop control $(v_t)_{s \leq t \leq T} \triangleq (u(t, X_{\cdot}^{s,x,\alpha,u}))_{s \leq t \leq T} \in \mathcal{U}_s$ which shows the inclusion above. 
\ep
}
\end{Remark}

We have the following well-posedness result for equation \eqref{SDE} when $u$ is an elementary feedback strategy (so that $u_r$ stands for $u(r,X_\cdot,I_{\cdot^-})$), where the only difference with Proposition \ref{P:X} is that now the solution is adapted to the smaller filtration $\F^{W,s}$, since $\F^s$ plays no role when $u\in\Uc_s^E$.

\begin{Proposition}
\label{P:X_Elem}
Let Assumption {\bf (H1)} hold. For any $(s,x,i)\in[0,T]\times\R^d\times\I_m$, $\alpha\in\Ac_s$, $u\in\Uc_s^E$, there exists a unique (up to indistinguishability) $\F^{W,s}$-adapted process $(X^{s,x,i;\alpha,u},I^{s,x,i;\alpha,u})=(X_t^{s,x,i;\alpha,u},I_t^{s,x,i;\alpha,u})_{s\leq t\leq T}$ to equation \eqref{SDE}, such that every path of $(X_\cdot^{s,x,i;\alpha,u},I_{\cdot^-}^{s,x,i;\alpha,u})$ belongs to $C([s,T];\R^d)\times\mathscr{L}([s,T];\I_m)$. Moreover, for any $q\geq1$ there exists a positive constant $C_{q,T}$, depending only on $q,T,M_1$ (independent of $s,x,i,\alpha,u$), such that
\beq
\label{EstimateXbis}
\E\Big[\sup_{s\leq t\leq T}|X_t^{s,x,i;\alpha,u}|^q\Big] &\leq & C_{q,T}(1+|x|^q).
\enq
\end{Proposition}
\textbf{Proof.}
The proof can be done along the lines of the proof of Proposition \ref{P:X}. We simply notice that in Proposition \ref{P:X} we used the following result: if $u\in\Uc_s$ and $I=(I_t)_{s\leq t\leq\tau}$ is known up to a certain $\F^s$-stopping time $\tau$, then there exists a unique (up to indistinguishability) $\F^s$-adapted solution $X=(X_t)_{s\leq t\leq\tau}$ to the equation
\beq
\label{X_ElementaryStrategy}
X_t &=& x + \int_s^t b(X_r,I_r,u_r)dr + \int_s^t\sigma(X_r,I_r,u_r)dW_r, \qquad s\leq t\leq\tau,
\enq
such that every path of $X$ belongs to $C([s,T];\R^d)$. The validity of this result is well-known under {\bf (H1)}. On the other hand, it is not immediately clear when $u\in\Uc_s^E$ is an elementary feedback strategy. However, the result is still valid and follows from Proposition 2.4 in \cite{sirbu14a}, see also Theorem 2.2 in \cite{sirbu14b}. Moreover, when $u\in\Uc_s^E$ it turns out that the process $X$ is adapted to the smaller filtration $\F^{W,s}$. Finally, under Assumption {\bf (H1)}, estimate \eqref{EstimateXbis} is well-known, see, e.g., Theorem 1.3.15 in \cite{pham09}.
\ep

\vspace{3mm}

We can finally introduce the value function for the robust switching control problem where nature adopts the elementary feedback strategies:
\begin{align*}
\overline V(s,x,i) &:= \sup_{\alpha\in\Ac_s}\!\inf_{u\in\Uc_s^E}\!\E\bigg[\int_s^T \!\!f(X_t,I_t,u_t')dt + g(X_T,I_T) - \sum_{n\in\N}c(X_{\tau_n'},I_{(\tau_n')^-},I_{\tau_n'})1_{\{s\leq\tau_n'<T\}}\bigg],
\end{align*}
for every $(s,x,i)\in[0,T]\times\R^d\times\I_m$, with the shorthands $X=X^{s,x,i;\alpha,u}$, $I=I^{s,x,i;\alpha,u}$, $\tau_n'=\tau_n(X_\cdot,I_{\cdot^-})$, and $u_t'=u(t,X_\cdot,I_{\cdot^-})$. 
This auxiliary formulation of robust switching problem where both players use feedback strategies (or closed-loop controls) is the same as the one used in  \cite{sirbu14a}.
Notice that $u'\in\Uc_s$ and we have
\beqs
\overline V(s,x,i) &:=& \sup_{\alpha\in\Ac_s}\inf_{u\in\Uc_s^E} J(s,x,i;\alpha,u'), \qquad \forall\,(s,x,i)\in[0,T]\times\R^d\times\I_m.
\enqs
In particular, $V(s,x,i)\leq\overline V(s,x,i)$, for any $(s,x,i)\in[0,T]\times\R^d\times\I_m$. Moreover, proceeding as in Remark \ref{R:V_polynomial_growth}, we can show that $\overline V$ satisfies a polynomial growth condition in $x$: $|\overline V(s,x,i)|\leq C(1+|x|^p)<\infty$, for some positive constant $C$, depending only on $T,M_1,M_2$, and with the same $p$ as in Assumption {\bf (H2)}(iii).

\subsection{Concatenation of feedback strategies}

In the present section, we need to introduce the concept of feedback control starting at a certain stopping rule $\tau$ and to define the notion of concatenation at $\tau$ of two feedback controls, which will be crucial in the development of the stochastic Perron method.

\begin{Definition}[Feedback switching strategies starting strictly later than $\tau$]
\label{D:switching_at_tau}
Fix $s$ in $[0,T]$ and $\tau\in\Tc^s$. We say that the double sequence $\alpha=(\tau_n,\iota_n)_{n\in\N}$ is a feedback switching strategy starting strictly later than $\tau$ if $\alpha\in\Ac_s$, with $\tau\leq\tau_0$ and $\tau<\tau_0$ on the set $\{\tau<T\}$. $\Ac_{s,\tau^+}$ denotes the family of all feedback switching strategies  for the controller, given the initial deterministic time $s$ and starting strictly later than $\tau$. When $\tau\equiv s$, we simply write $\Ac_{s^+}$ instead of $\Ac_{s,s^+}$.
\end{Definition}

Following \cite{sirbu14a}, Definition 2.7, and recalling Remark \ref{R:Strategy_tau^+}, we now define the elementary feedback strategies  starting at some stopping rule $\tau$.

\begin{Definition}[Elementary feedback strategies  starting at $\tau$]
\label{D:ElemStrategy_at_tau}
Fix $s\in[0,T]$ and $\tau\in\Tc^s$. We say that $u$ is an elementary feedback strategy starting at $\tau$ if:
\begin{itemize}
\item $\tau_k\in\Tc^s$, for any $k=1,\ldots,n$, and
\[
\tau \ =: \ \tau_0 \ \leq \ \cdots \ \leq \ \tau_k \ \leq \ \cdots \ \leq \ \tau_n \ = \ T.
\]
\item $\xi_k\colon C([s,T];\R^d)\times\mathscr{L}([s,T];\I_m)\rightarrow U$ is $\Bc_{\tau_{k-1}^+}^s$-measurable, for any $k=1,\ldots,n$.
\end{itemize}
The elementary feedback strategy 
\beqs
u\colon\big\{(t,y)\in[s,T]\times\big(C([s,T];\R^d)\times\mathscr{L}([s,T];\I_m)\big)\colon\tau(y)\leq t\leq T\big\} &\longrightarrow & U
\enqs
is given by
\[
u(t,y) \ := \ \xi_1(y)1_{\{t=\tau(y)\}} + \sum_{k=1}^n \xi_k(y)1_{\{\tau_{k-1}(y)<t\leq\tau_k(y)\}}.
\]
$\Uc_{s,\tau}^E$ denotes the family of all elementary feedback strategy, given the initial deterministic time $s$ and starting at $\tau$.
\end{Definition}

Notice that, when $\tau=s$ in Definition \ref{D:ElemStrategy_at_tau}, the set $\Uc_s^E$ is just $\Uc_s^E$.

\begin{Remark}
\label{R:strictly}
{\rm
Definition \ref{D:switching_at_tau} is inspired by Definition 2.7 in \cite{sirbu14a}, with in addition the condition ``\,$\tau<\tau_0$ on the set $\{\tau<T\}$'', which justifies the presence of the adverb \emph{strictly} in the name. Indeed, our aim is to define the set $\Ac_{s,\tau^+}$ in such a way that when we concatenate two feedback switching strategies $\alpha\in\Ac_s$ and $\tilde\alpha\in\Ac_{s,\tau^+}$ at a stopping rule $\tau\in\Tc^s$ (see Proposition \ref{P:Concatenate} below) then $\alpha\otimes_\tau\tilde\alpha$ coincides with $\alpha$ at time $\tau$ (this property plays an important role in the sequel, e.g., in the proof of Theorem \ref{T:StochPerron}). On the other hand, when we concatenate two elementary feedback strategies $u\in\Uc_s^E$ and $\tilde u\in\Uc_{s,\tau}^E$, then $u\otimes_\tau\tilde u$ coincides with $u$ at time $\tau$, simply adopting the same definition for $\Uc_{s,\tau}^E$ as in \cite{sirbu14a} combined with Remark \ref{R:Strategy_tau^+}.
\ep
}
\end{Remark}

As in \cite{sirbu14a}, Lemma 2.8 and Proposition 2.9, we have the two following results, whose simple proof is only sketched for Lemma \ref{L:Partition} and omitted for Proposition \ref{P:Concatenate}.

\begin{Lemma}
\label{L:Partition}
Fix $s\in[0,T]$, $\tau\in\Tc^s$, $\alpha^1=(\tau_n^1,\iota_n^1)_{n\in\N},\alpha^2=(\tau_n^2,\iota_n^2)_{n\in\N}\in\Ac_{s,\tau^+}$, $u^1,u^2\in\Uc_{s,\tau}^E$, and $B\in\Bc_{\tau^+}^s$.
\begin{itemize}
\item The double sequence $\alpha=(\tau_n,\iota_n)_{n\in\N}$ given by
\beqs
\big(\tau_n(y),\iota_n(y)\big) &=& \big(\tau_n^1(y),\iota_n^1(y)\big)1_{\{y\in B\}} + \big(\tau_n^2(y),\iota_n^2(y)\big)1_{\{y\in B^c\}}
\enqs
is in $\Ac_{s,\tau^+}$.
\item The map
\beqs
u\colon\big\{(t,y)\in[s,T]\times\big(C([s,T];\R^d)\times\mathscr{L}([s,T];\I_m)\big)\colon\tau(y)\leq t\leq T\big\} &\longrightarrow & U
\enqs
given by
\beqs
u(t,y) &=& u^1(t,y)1_{\{y\in B\}} + u^2(t,y)1_{\{y\in B^c\}}
\enqs
is in $\Uc_{s,\tau}^E$.
\end{itemize}
\end{Lemma}
\textbf{Proof.}
We only prove the first item, where we focus on the two main points. In particular, the proof that $\tau_n\in\Tc^s$ and $\iota_n\in\Bc_{\tau_n}^s$ is based on the observation that $B\in\Bc_{\tau^+}^s\subset\Bc_{\tau_n^1}^s,\Bc_{\tau_n^2}^s$, for any $n\in\N$, which is a consequence of the property: $\tau<\tau_0^1,\tau_0^2$ on the set $\{\tau<T\}$. The other non-trivial part is the proof that $\alpha$ satisfies property (ii) of Remark \ref{R:Property_tau_n}. To prove it, consider $(y_n)_{n\in\N}\in C([s,T];\R^d)\times\mathscr{L}([s,T];\I_m)$, with $y_n(t)=y_{n+1}(t)$, $t\in[s,\tau_n(y_n)]$. Since $\tau_0\leq\tau_n$ for any $n\in\N$, we have
\beqs
y_0(t) &=& y_n(t), \qquad \forall\,t\in[s,\tau_0(y_0)],\,n\in\N.
\enqs
As $\tau<\tau_0$ on the set $\{\tau<T\}$, it follows that
\beq
\label{E:y_0=y_n}
y_0(t^+) &=& y_n(t^+), \qquad \forall\,t\in[s,\tau(y_0)],\,n\in\N.
\enq
In particular $y_0(\tau(y_0)^+)=y_n(\tau(y_0)^+)$. Moreover, from Lemma \ref{L:tau&iota1} we get $\tau(y_0)=\tau(y_n)$, so that $y_0(\tau(y_0)^+)=y_n(\tau(y_n)^+)$. Therefore, $y_0\in B$ if and only if $y_n\in B$, for any $n\in\N$. In conclusion, property (ii) of Remark \ref{R:Property_tau_n} for $(\tau_n)_{n\in\N}$ follows from the definitions of $(\tau_n^1)_{n\in\N}$ and $(\tau_n^2)_{n\in\N}$.
\ep

\begin{Proposition}[Concatenation]
\label{P:Concatenate}
Fix $s\in[0,T]$, $\tau,\rho\in\Tc^s$ with $\tau\leq\rho\leq T$, $\tilde\alpha=(\tilde\tau_n,\tilde\iota_n)_{n\in\N}\in\Ac_{s,\rho^+}$, $\tilde u\in\Uc_{s,\rho}^E$. Then
\begin{itemize}
\item for each $\alpha=(\tau_n,\iota_n)_{n\in\N}\in\Ac_s$ $($resp. $\alpha=(\tau_n,\iota_n)_{n\in\N}\in\Ac_{s,\tau^+}$$)$, the double sequence $\alpha\otimes_\rho\tilde\alpha=(\tau_n^{\otimes_\rho},\iota_n^{\otimes_\rho})_{n\in\N}$ given by
\beqs
\big(\tau_n^{\otimes_\rho}(y),\iota_n^{\otimes_\rho}(y)\big) &=& \big(\tau_n(y),\iota_n(y)\big)1_{\{\tau_n(y)\leq\rho(y)\}} + \big(\tilde\tau_n(y),\tilde\iota_n(y)\big)1_{\{\tau_n(y)>\rho(y)\}}
\enqs
is in $\Ac_s$ $($resp. $\Ac_{s,\tau^+}$$)$;
\item for each $u\in\Uc_{s,\tau}^E$, the map
\beqs
u\otimes_\rho\tilde u\colon\big\{(t,y)\in[s,T]\times\big(C([s,T];\R^d)\times\mathscr{L}([s,T];\I_m)\big)\colon\tau(y)\leq t\leq T\big\} &\longrightarrow & U
\enqs
given by
\beqs
(u\otimes_\rho\tilde u)(t,y) &=& u(t,y)1_{\{\tau(y)\leq t\leq\rho(y)\}} + \tilde u(t,y)1_{\{\rho(y)<t\leq T\}}
\enqs
is in $\Uc_{s,\tau}^E$.
\end{itemize}
\end{Proposition}

\subsection{Definitions of $\Vc^-$, $\Vc^+$ and their properties}

We can now provide the definitions of the classes of functions $\Vc^-$ and $\Vc^+$, which are the cornerstones of the stochastic Perron method. Their elements are known in the literature on stochastic Perron method as stochastic subsolutions ($\Vc^-$) and stochastic supersolutions ($\Vc^+$), see e.g. \cite{BS13}. 

\begin{Definition}
\label{D:StochSubSol}
$\Vc^-$ is the set of functions $v\colon[0,T]\times\R^d\times\I_m\rightarrow\R$ which have the following properties:
\begin{itemize}
\item $v$ is continuous and satisfies the terminal condition $v(T,x,i)\leq g(x,i)$, $(x,i)\in\R^d\times\I_m$, together with the polynomial growth condition
\beqs
\sup_{(s,x,i)\in[0,T]\times\R^d\times\I_m}\frac{|v(s,x,i)|}{1+|x|^q} &<& \infty,
\enqs
for some $q\geq 1$.
\item For any $s\in[0,T]$ and $\tau,\rho\in\Tc^s$ with $\tau\leq\rho\leq T$, there exists $\tilde\alpha=(\tilde\tau_n,\tilde\iota_n)_{n\in\N}\in\Ac_{s,\tau^+}$ (possibly depending on $s,\tau,\rho$) such that, for any $\alpha=(\tau_n,\iota_n)_{n\in\N}\in\Ac_s$, $u\in\Uc_s$, and $(x,i)\in\R^d\times\I_m$, we have
\beqs
v(\tau',X_{\tau'},I_{\tau'}) &\leq & \E\bigg[\int_{\tau'}^{\rho'} f(X_t,I_t,u_t)dt + v(\rho',X_{\rho'},I_{\rho'}) \\
& & \quad - \; \sum_{n\in\N}c(X_{\tilde\tau_n'},I_{(\tilde\tau_n')^-},I_{\tilde\tau_n'})1_{\{\tau'\leq\tilde\tau_n'<\rho'\}}\bigg|\Fc_{\tau'}^s\bigg], \qquad \P\text{-a.s.}
\enqs
with the shorthands $X=X^{s,x,i;\alpha\otimes_\tau\tilde\alpha,u}$, $I=I^{s,x,i;\alpha\otimes_\tau\tilde\alpha,u}$, $\tau'=\tau(X_\cdot,I_{\cdot^-})$, $\rho'=\rho(X_\cdot,I_{\cdot^-})$, and $\tilde\tau_n'=\tilde\tau_n(X_\cdot,I_{\cdot^-})$.
\end{itemize}
\end{Definition}

\begin{Definition}
\label{D:StochSuperSol}
$\Vc^+$ is the set of functions $v\colon[0,T]\times\R^d\times\I_m\rightarrow\R$ which have the following properties:
\begin{itemize}
\item $v$ is continuous and satisfies the terminal condition $v(T,x,i)\geq g(x,i)$, $(x,i)\in\R^d\times\I_m$, together with the polynomial growth condition
\beqs
\sup_{(s,x,i)\in[0,T]\times\R^d\times\I_m}\frac{|v(s,x,i)|}{1+|x|^q} &<& \infty,
\enqs
for some $q\geq1$.
\item For any $s\in[0,T]$, $\tau\in\Tc^s$, and $\alpha=(\tau_n,\iota_n)_{n\in\N}\in\Ac_s$, there exists $\tilde u\in\Uc_{s,\tau}^E$ (possibly depending on $s,\tau,\alpha$) such that, for any $u\in\Uc_s^E$, $(x,i)\in\R^d\times\I_m$, and $\rho\in\Tc^s$, with $\tau\leq\rho\leq T$, we have
\beqs
v(\tau',X_{\tau'},I_{\tau'}) &\geq & \E\bigg[\int_{\tau'}^{\rho'} f(X_t,I_t,\tilde u_t)dt + v(\rho',X_{\rho'},I_{\rho'}) \\
& & \quad - \; \sum_{n\in\N}c(X_{\tau_n'},I_{(\tau_n')^-},I_{\tau_n'})1_{\{\tau'\leq\tau_n'<\rho'\}}\bigg|\Fc_{\tau'}^s\bigg], \qquad \P\text{-a.s.}
\enqs
with the shorthands $X=X^{s,x,i;\alpha,u\otimes_\tau\tilde u}$, $I=I^{s,x,i;\alpha,u\otimes_\tau\tilde u}$, $\tau'=\tau(X_\cdot,I_{\cdot^-})$, $\rho'=\rho(X_\cdot,I_{\cdot^-})$, $\tau_n'=\tau_n(X_\cdot,I_{\cdot^-})$, and $\tilde u_t=\tilde u(t,X_\cdot,I_{\cdot^-})$.
\end{itemize}
\end{Definition}

\begin{Remark}
{\rm
The definitions of $\Vc^-$ and $\Vc^+$ are inspired by \cite{sirbu14a}, Definitions 3.1-3.2-3.3, but for the fact that in Definition \ref{D:StochSubSol} above we fix $\rho$ before $\tilde\alpha$, so that $\tilde\alpha$ can depend on $\rho$. This greater freedom in the choice of $\tilde\alpha$ turns out to be fundamental in the implementation of the stochastic Perron method, Theorem \ref{T:StochPerron}, and it is due to the condition ``\,$\tau<\tau_0$ on the set $\{\tau<T\}$'' in the definition of $\Ac_{s,\tau^+}$, already discussed in Remark \ref{R:strictly}. Indeed, using the set $\Ac_{s,\tau^+}$, the existence of an ``optimal'' feedback switching strategy  $\tilde\alpha=(\tilde\tau_n,\tilde\iota_n)_{n\in\N}\in\Ac_{s,\tau^+}$, which works for every $\rho\in\Tc^s$, with $\tau\leq\rho\leq T$, is not guaranteed. For example, it could happen that every ``optimal'' feedback switching strategy which works for all $\rho$ has to satisfy $\tilde\tau_0=\tau$, therefore it can not belong to $\Ac_{s,\tau^+}$. To avoid this problem, firstly we fix $\rho$, then we choose an ``optimal'' $\tilde\alpha\in\Ac_{s,\tau^+}$. Another possibility would be to look for an ``$\eps$-optimal'' $\tilde\alpha\in\Ac_{s,\tau^+}$ which works for every $\rho$.
\ep
}
\end{Remark}

We first notice that, as stated below, the two sets $\Vc^-$ and $\Vc^+$ are not empty, moreover every $v\in\Vc^-$ (resp. $v\in\Vc^+$) satisfies the sub-dynamic (resp. super-dynamic) programming principle, also known as suboptimality (resp. superoptimality) principle, see \cite{swiech96}.

\begin{Lemma}
\label{L:V^-nonempty_HalfDPP}
Let Assumptions {\bf (H1)} and {\bf (H2)} hold.
\begin{itemize}
\item[\textup{(i)}] $\Vc^-\neq\emptyset$ and $\Vc^+\neq\emptyset$.
\item[\textup{(ii)}] Every $v\in\Vc^-$ satisfies the sub-dynamic programming principle: for any $(s,x,i)\in[0,T]\times\R^d\times\I_m$ and $\rho\in\Tc^s$,
\beq
\label{DPP^-}
v(s,x,i) &\leq & \sup_{\alpha\in\Ac_{s^+}}\inf_{u\in\Uc_s} \E\bigg[\int_s^{\rho'} f(X_t,I_t,u_t)dt + v(\rho',X_{\rho'},I_{\rho'}) \\
& & \hspace{2.5cm} - \; \sum_{n\in\N}c(X_{\tau_n'},I_{(\tau_n')^-},I_{\tau_n'})1_{\{s\leq\tau_n'<\rho'\}}\bigg], \notag
\enq
with the shorthands $X=X^{s,x,i;\alpha,u}$, $I=I^{s,x,i;\alpha,u}$, $\rho'=\rho(X_\cdot,I_{\cdot^-})$, and $\tau_n'=\tau_n(X_\cdot,I_{\cdot^-})$.
\item[\textup{(iii)}] Every $v\in\Vc^+$ satisfies the super-dynamic programming principle: for any $(s,x,i)\in[0,T]\times\R^d\times\I_m$ and $\rho\in\Tc^s$,
\beq
\label{DPP^+}
v(s,x,i) &\geq & \sup_{\alpha\in\Ac_{s^+}}\inf_{u\in\Uc_s^E} \E\bigg[\int_s^{\rho'} f(X_t,I_t,u_t)dt + v(\rho',X_{\rho'},I_{\rho'}) \\
& & \hspace{2.5cm} - \; \sum_{n\in\N}c(X_{\tau_n'},I_{(\tau_n')^-},I_{\tau_n'})1_{\{s\leq\tau_n'<\rho'\}}\bigg], \notag
\enq
with the shorthands $X=X^{s,x,i;\alpha,u}$, $I=I^{s,x,i;\alpha,u}$, $\rho'=\rho(X_\cdot,I_{\cdot^-})$, $\tau_n'=\tau_n(X_\cdot,I_{\cdot^-})$, and $u_t=u(t,X_\cdot,I_{\cdot^-})$.
\end{itemize}
\end{Lemma}
\textbf{Proof.}
We begin proving that $\Vc^-\neq\emptyset$. Let us consider the function $v\colon[0,T]\times\R^d\times\I_m\rightarrow\R$ given by
\beq
\label{v^*}
v(s,x,i) &:=& -Ce^{\lambda(T-s)}(1 + |x|^q), \qquad \forall\,(s,x,i)\in[0,T]\times\R^d\times\I_m,
\enq
where $q=\max\{4,p\}$, with $p$ as in Assumption {\bf (H2)}(iii), and $C,\lambda$ are positive constants to be determined later. Set $h(x)=|x|^q$. Notice that $h\in C^2(\R^d)$ and there exists a positive constant $M_h$ (depending only on $q$) such that $|D_x h(x)|\leq M_h|x|^{q-1}$ and $D_x^2 h(x)\leq M_h|x|^{q-2}$, $\forall\,x\in\R^d$.

From the polynomial growth condition of $g$ in Assumption {\bf (H2)}(iii), we see that $v(T,x,i)\leq g(x,i)$ if we choose $C$ large enough.

Now, we choose $\lambda$ opportunely. Fix $s\in[0,T]$ and $\tau,\rho\in\Tc^s$ with $\tau\leq\rho\leq T$. We choose $\tilde\alpha=(\tilde\tau_n,\tilde\iota_n)_{n\in\N}\in\Ac_{s,\tau^+}$ as follows: for any $n\in\N$, $\tilde\tau_n\equiv T$ and $\tilde\iota_n\equiv\underline i$, for some fixed $\underline i\in\I_m$. Let $\alpha=(\tau_n,\iota_n)_{n\in\N}\in\Ac_s$, $u\in\Uc_s$, and $(x,i)\in\R^d\times\I_m$. Set $X=X^{s,x,i;\alpha\otimes_\tau\tilde\alpha,u}$, $I=I^{s,x,i;\alpha\otimes_\tau\tilde\alpha,u}$, $\tau'=\tau(X,I)$, and $\rho'=\rho(X,I)$. Then, noting that $v(r,X_r,I_r)$ is constant with respect to $I_r$, and applying It\^o's formula to $\int_{\tau'}^r f(X_t,I_t,u_t)dt + v(r,X_r,I_r)$ between $\tau'$ and $\rho'$, we obtain
\beq
& & \int_{\tau'}^{\rho'} f(X_t,I_t,u_t)dt + v(\rho',X_{\rho'},I_{\rho'}) \notag \\
&=& \int_{\tau'}^{\rho'} f(X_t,I_t,u_t)dt + v(\tau',X_{\tau'},I_{\tau'}) - C\int_{\tau'}^{\rho'} e^{\lambda(T-t)} D_x h(X_t).b(X_t,I_t,u_t)dt \notag \\
& & - \; C\int_{\tau'}^{\rho'} e^{\lambda(T-t)} (D_x h(X_t))\trans \sigma(X_t,I_t,u_t) dW_t + \lambda C\int_{\tau'}^{\rho'}e^{\lambda(T-t)}(1+h(X_t))dt \notag \\
& & - \; \frac{1}{2}C\int_{\tau'}^{\rho'} e^{\lambda(T-t)}\text{tr}\big[\sigma\sigma\trans(X_t,I_t,u_t)D_x^2 h(X_t)\big]dt. \label{Proof1}
\enq
Consider the $\F^s$-local martingale $M_r=\int_s^r 1_{[\tau',T]}(t)e^{\lambda(T-t)} (D_x h(X_t))\trans \sigma(X_t,I_t,u_t) dW_t$, $r\in[s,T]$. In order to prove that $M$ is a true martingale, we show that $\E[\sup_{s\leq r\leq T}|M_r|]<\infty$. From Burkholder-Davis-Gundy's inequality, we see that it is enough to prove $\E[\sqrt{\langle M\rangle_T}]<\infty$, namely
\[
\E\bigg[\sqrt{\int_{\tau'}^T e^{2\lambda(T-t)} |D_x h(X_t)|^2 \|\sigma(X_t,I_t,u_t)\|^2 dt}\bigg] \ < \ \infty.
\]
This latter inequality holds since $|D_xh(x)|\leq M_h|x|^{q-1}$, $\|\sigma(x,i,u)\|\leq M_1(1+|x|)$ (see Remark \ref{R:LinearGrowth_b_sigma}), and $X$ satisfies estimate \eqref{EstimateX}. From the martingale property of $M$ and Doob's optional sampling theorem, we have in particular
\[
\E\bigg[\int_{\tau'}^{\rho'} e^{\lambda(T-t)} (D_x h(X_t))\trans \sigma(X_t,I_t,u_t) dW_t\bigg|\Fc_{\tau'}^s\bigg] \ = \ \E\big[M_{\rho'}\big|\Fc_{\tau'}^s\big] \ = \ 0.
\]
Therefore, taking the conditional expectation with respect to $\Fc_{\tau'}^s$ in \eqref{Proof1}, using the linear growth conditions of $b,\sigma,f$, and the estimates on $D_x h(x)$ and $D_x^2 h(x)$, we find
\beqs
& & \E\bigg[\int_{\tau'}^{\rho'} f(X_t,I_t,u_t)dt + v(\rho',X_{\rho'},I_{\rho'})\bigg|\Fc_{\tau'}^s\bigg] \\
&\geq & v(\tau',X_{\tau'},I_{\tau'}) + \E\bigg[- M_2\int_{\tau'}^{\rho'} (1+|X_t|^p)dt - C M_h M_1\int_{\tau'}^{\rho'} e^{\lambda(T-t)} |X_t|^{q-1}(1+|X_t|)dt \\
& & + \lambda C\int_{\tau'}^{\rho'}e^{\lambda(T-t)}(1+|X_t|^q)dt - \frac{1}{2}C M_h M_1^2\int_{\tau'}^{\rho'} e^{\lambda(T-t)}|X_t|^{q-2}(1+|X_t|)^2 dt \bigg|\Fc_{\tau'}^s\bigg].
\enqs
We see that there exists a positive constant $\bar C$ (depending only on $C,M_h,M_1,M_2$) such that
\beqs
& & \E\bigg[\int_{\tau'}^{\rho'} f(X_t,I_t,u_t)dt + v(\rho',X_{\rho'},I_{\rho'})\bigg|\Fc_{\tau'}^s\bigg] \;\, \geq \;\, v(\tau',X_{\tau'},I_{\tau'}) \\
& & + \; (\lambda C - \bar C)\,\E\bigg[\int_{\tau'}^{\rho'} e^{\lambda(T-t)}(1+|X_t|^q)dt\bigg|\Fc_{\tau'}^s\bigg].
\enqs
Now, we choose $\lambda\geq0$ such that $\lambda C-\bar C\geq0$. Then, we have
\beqs
\E\bigg[\int_{\tau'}^{\rho'} f(X_t,I_t,u_t)dt + v(\rho',X_{\rho'},I_{\rho'})\bigg|\Fc_{\tau'}^s\bigg] &\geq & v(\tau',X_{\tau'},I_{\tau'}).
\enqs
From the definition of $\tilde\alpha$, we see that $\sum_{n\in\N}c(X_{\tilde\tau_n'},I_{(\tilde\tau_n')^-},I_{\tilde\tau_n'})1_{\{\tau'\leq\tilde\tau_n'<\rho'\}}=0$. Therefore, it follows that $v\in\Vc^-$. In a similar way we can prove that $-v\in\Vc^+$, so that $\Vc^+\neq\emptyset$.

Concerning (ii), let $v\in\Vc^-$ and fix $s\in[0,T]$, $\tau,\rho\in\Tc^s$, with $s\equiv\tau\leq\rho\leq T$. From the second item of the definition of $\Vc^-$, there exists $\tilde\alpha=(\tilde\tau_n,\tilde\iota_n)_{n\in\N}\in\Ac_{s^+}$ such that, for any $u\in\Uc_s$ and $(x,i)\in\R^d\times\I_m$ (we choose $\alpha=(\tau_n,\iota_n)_{n\in\N}\in\Ac_s$ with $\tau_n\equiv T$ and $\iota_n\equiv i$, for any $n\in\N$; with this choice we have $(X^{s,x,i;\alpha\otimes_s\tilde\alpha,u},I^{s,x,i;\alpha\otimes_s\tilde\alpha,u})=(X^{s,x,i;\tilde\alpha,u},I^{s,x,i;\tilde\alpha,u})$; in particular, $I_s^{s,x,i;\alpha\otimes_s\tilde\alpha,u}=i$), we find
\beq
\label{HalfDPP_V^-_proof}
v(s,x,i) &\leq & \E\bigg[\int_s^{\rho'} f(X_t,I_t,u_t)dt + v(\rho',X_{\rho'},I_{\rho'}) \notag \\
& & \quad - \; \sum_{n\in\N}c(X_{\tilde\tau_n'},I_{(\tilde\tau_n')^-},I_{\tilde\tau_n'})1_{\{s\leq\tilde\tau_n'<\rho'\}}\bigg|\Fc_s^s\bigg], \qquad \P\text{-a.s.}
\enq
with the shorthands $X=X^{s,x,i;\tilde\alpha,u}$, $I=I^{s,x,i;\tilde\alpha,u}$, $\rho'=\rho(X_\cdot,I_{\cdot^-})$, and $\tilde\tau_n'=\tilde\tau_n(X_\cdot,I_{\cdot^-})$. Taking the expectation in \eqref{HalfDPP_V^-_proof} and the infimum with respect to $u\in\Uc_s$, we get
\beqs
v(s,x,i) &\leq & \inf_{u\in\Uc_s}\E\bigg[\int_s^{\rho'} f(X_t,I_t,u_t)dt + v(\rho',X_{\rho'},I_{\rho'}) \\
& & \hspace{2cm} - \; \sum_{n\in\N}c(X_{\tilde\tau_n'},I_{(\tilde\tau_n')^-},I_{\tilde\tau_n'})1_{\{s\leq\tilde\tau_n'<\rho'\}}\bigg] \\
&\leq & \sup_{\alpha\in\Ac_{s^+}}\inf_{u\in\Uc_s}\E\bigg[\int_s^{\rho'} f(X_t,I_t,u_t)dt + v(\rho',X_{\rho'},I_{\rho'}) \\
& & \hspace{2cm} - \; \sum_{n\in\N}c(X_{\tau_n'},I_{(\tau_n')^-},I_{\tau_n'})1_{\{s\leq\tau_n'<\rho'\}}\bigg].
\enqs
In a similar way we can prove statement (iii).
\ep

\vspace{3mm}

As stated below, every $v\in\Vc^-$ is less than every $v\in\Vc^+$, while the value functions $V$ and $\overline V$ are squeezed between them.

\begin{Lemma}
Let Assumptions {\bf (H1)} and {\bf (H2)} hold.
\begin{itemize}
\item[\textup{(i)}] $\sup_{v\in\Vc^-}v$ $=:$ $v^-$ $\leq$ $V$ $\leq$ $\overline V$ $\leq$ $v^+$ $:=$ $\inf_{v\in\Vc^+}v$.
\item[\textup{(ii)}] $v^-$ is lsc and satisfies the polynomial growth condition
\beq
\label{v^-GrowthCond}
\sup_{(s,x,i)\in[0,T]\times\R^d\times\I_m}\frac{|v^-(s,x,i)|}{1+|x|^q} &<& \infty,
\enq
for some $q\geq1$.
\item[\textup{(iii)}] $v^+$ is usc and satisfies the polynomial growth condition
\beqs
\sup_{(s,x,i)\in[0,T]\times\R^d\times\I_m}\frac{|v^+(s,x,i)|}{1+|x|^q} &<& \infty,
\enqs
for some $q\geq1$.
\end{itemize}
\end{Lemma}
\textbf{Proof.}
Concerning (i), to obtain the inequality $v\leq V$ for all $v\in\Vc^-$ (resp. $\overline V\leq v$ for all $v\in\Vc^+$) we take $\rho\equiv T$ in the sub-dynamic programming principle \eqref{DPP^-} (resp. super-dynamic programming principle \eqref{DPP^+}) and we use the inequality $v(T,x,i)\leq g(x,i)$ (resp. $v(T,x,i)\geq g(x,i)$), for all $(x,i)\in\R^d\times\I_m$. Regarding (ii), we notice that $v^-$ is lsc since it is the supremum of a family of lsc (actually, continuous) functions. Moreover, let $\underline v\in\Vc^-$ and $\bar v\in\Vc^+$. From (i) it follows that $\underline v\leq v^-\leq \bar v$, and from the polynomial growth condition of $\underline v,\bar v$ we see that $v^-$ satisfies the polynomial growth condition \eqref{v^-GrowthCond}. Statement (iii) can be proved in a similar way.
\ep

\vspace{3mm}

We can now state our main result.

\begin{Theorem}[Stochastic Perron's method]
\label{T:StochPerron}
Let Assumptions {\bf (H1)} and {\bf (H2)} hold. Then, $v^-$ is a viscosity supersolution to equation \eqref{HJB} and $v^+$ is a viscosity subsolution to equation \eqref{HJB}.
\end{Theorem}

In order to prove Theorem \ref{T:StochPerron}, we need the following two lemmata. In particular, Lemma \ref{L:Stability} states that $\Vc^-$ (resp. $\Vc^+$) is stable by supremum (resp. infimum), which gives the existence of a monotone approximating sequence for $v^-$ (resp. $v^+$) in Lemma \ref{L:ApproximatingSequence}.

\begin{Lemma}
\label{L:Stability}
Let Assumptions {\bf (H1)} and {\bf (H2)} hold.
\begin{itemize}
\item[\textup{(i)}] If $v^1,v^2\in\Vc^-$ then $v:=v^1\vee v^2\in\Vc^-$.
\item[\textup{(ii)}] If $v^1,v^2\in\Vc^+$ then $v:=v^1\wedge v^2\in\Vc^+$.
\end{itemize}
\end{Lemma}
\textbf{Proof.}
Let us prove (i). As the first item in Definition \ref{D:StochSubSol} clearly holds, we prove that $v$ satisfies the second item. To this end, fix $s\in[0,T]$ and $\tau,\rho\in\Tc^s$ with $\tau\leq\rho\leq T$. Let $\tilde\alpha^1=(\tilde\tau_n^1,\tilde\iota_n^1)_{n\in\N},\tilde\alpha^2=(\tilde\tau_n^2,\tilde\iota_n^2)_{n\in\N}\in\Ac_{s,\tau^+}$ be the two feedback switching controls, starting strictly later than $\tau$, corresponding to $v^1$ and $v^2$. Now, consider the set $B:=\{(v^1-v^2)(\tau(y),y(\tau(y)^+))\geq0\}\in\Bc_{\tau^+}^s$ and define the double sequence $\tilde\alpha=(\tilde\tau_n,\tilde\iota_n)_{n\in\N}$ as follows
\beqs
\big(\tilde\tau_n(y),\tilde\iota_n(y)\big) &:=& \big(\tilde\tau_n^1(y),\tilde\iota_n^1(y)\big) 1_{\{y\in B\}} + \big(\tilde\tau_n^2(y),\tilde\iota_n^2(y)\big) 1_{\{y\in B^c\}},
\enqs
for any $y\in C([s,T];\R^d)\times\mathscr{L}([s,T];\I_m)$, $n\in\N$. From Lemma \ref{L:Partition} it follows that $\tilde\alpha\in\Ac_{s,\tau^+}$. Now, we prove that $\tilde\alpha$ satisfies the condition in the second item of Definition \ref{D:StochSubSol}. Take $\alpha=(\tau_n,\iota_n)_{n\in\N}\in\Ac_s$, $u\in\Uc_s$, and $(x,i)\in\R^d\times\I_m$. We adopt the shorthands:
\[
\begin{array}{ccc}
X \ = \ X^{s,x,i;\alpha\otimes_\tau\tilde\alpha,u}, & \qquad X^1 \ = \ X^{s,x,i;\alpha\otimes_\tau\tilde\alpha^1,u}, & \qquad X^2 \ = \ X^{s,x,i;\alpha\otimes_\tau\tilde\alpha^2,u}, \\
I \ = \ I^{s,x,i;\alpha\otimes_\tau\tilde\alpha,u}, & \qquad I^1 \ = \ I^{s,x,i;\alpha\otimes_\tau\tilde\alpha^1,u}, & \qquad I^2 \ = \ I^{s,x,i;\alpha\otimes_\tau\tilde\alpha^2,u}.
\end{array}
\]
We also denote $\tau'=\tau(X_\cdot,I_{\cdot^-})$, $\rho'=\rho(X_\cdot,I_{\cdot^-})$, $\rho^{1,'}=\rho(X_\cdot^1,I_{\cdot^-}^1)$, $\rho^{2,'}=\rho(X_\cdot^2,I_{\cdot^-}^2)$, $\tilde\tau_n'=\tilde\tau_n(X_\cdot,I_{\cdot^-})$, $\tilde\tau_n^{1,'}=\tilde\tau_n^1(X_\cdot^1,I_{\cdot^-}^1)$, and $\tilde\tau_n^{2,'}=\tilde\tau_n^2(X_\cdot^2,I_{\cdot^-}^2)$. Notice that $(X_t,I_{t^-})=(X_t^1,I_{t^-}^1)=(X_t^2,I_{t^-}^2)$, $t\in[s,\tau']$. Therefore, from Lemma \ref{L:tau&iota1} we see that $\tau'=\tau(X_\cdot^1,I_{\cdot^-}^1)=\tau(X_\cdot^2,I_{\cdot^-}^2)$. Moreover, for any $t\in[\tau',T]$,
\beqs
(X_t,I_t) &=& (X_t^1,I_t^1) 1_{\{(v^1-v^2)(\tau',X_{\tau'},I_{\tau'})\geq0\}} + (X_t^2,I_t^2) 1_{\{(v^1-v^2)(\tau',X_{\tau'},I_{\tau'})<0\}}.
\enqs
As a consequence,
\beqs
\rho' &=& \rho^{1,'} 1_{\{(v^1-v^2)(\tau',X_{\tau'},I_{\tau'})\geq0\}} + \rho^{2,'} 1_{\{(v^1-v^2)(\tau',X_{\tau'},I_{\tau'})<0\}}, \\
\tilde\tau_n' &=& \tilde\tau_n^{1,'} 1_{\{(v^1-v^2)(\tau',X_{\tau'},I_{\tau'})\geq0\}} + \tilde\tau_n^{2,'} 1_{\{(v^1-v^2)(\tau',X_{\tau'},I_{\tau'})<0\}}.
\enqs
Therefore, from the previous identities and the properties of $v^1$, we obtain
\begin{align*}
&v^1(\tau',X_{\tau'},I_{\tau'})1_{\{(v^1-v^2)(\tau',X_{\tau'},I_{\tau'})\geq0\}} \ = \ v^1(\tau',X_{\tau'}^1,I_{\tau'}^1)1_{\{(v^1-v^2)(\tau',X_{\tau'},I_{\tau'})\geq0\}} \\
&\leq \ \E\bigg[\bigg(\int_{\tau'}^{\rho^{1,'}} f(X_t^1,I_t^1,u_t)dt + v^1(\rho^{1,'},X_{\rho^{1,'}}^1,I_{\rho^{1,'}}^1) \\
&\quad - \sum_{n\in\N}c(X_{\tilde\tau_n^{1,'}}^1,I_{(\tilde\tau_n^{1,'})^-}^1,I_{\tilde\tau_n^{1,'}}^1)1_{\{\tau'\leq\tilde\tau_n^{1,'}<\rho^{1,'}\}}\bigg)1_{\{(v^1-v^2)(\tau',X_{\tau'},I_{\tau'})\geq0\}}\bigg|\Fc_{\tau'}^s\bigg] \\
&\leq \ \E\bigg[\bigg(\int_{\tau'}^{\rho'} f(X_t,I_t,u_t)dt + v(\rho',X_{\rho'},I_{\rho'}) \\
&\quad - \sum_{n\in\N}c(X_{\tilde\tau_n'},I_{(\tilde\tau_n')^-},I_{\tilde\tau_n'})1_{\{\tau'\leq\tilde\tau_n'<\rho'\}}\bigg)1_{\{(v^1-v^2)(\tau',X_{\tau'},I_{\tau'})\geq0\}}\bigg|\Fc_{\tau'}^s\bigg].
\end{align*}
Concerning $v^2$, proceeding similarly we get
\begin{align*}
v^2(\tau',X_{\tau'},I_{\tau'})1_{\{(v^1-v^2)(\tau',X_{\tau'},I_{\tau'})<0\}} \ \leq \ \E\bigg[\bigg(\int_{\tau'}^{\rho'} f(X_t,I_t,u_t)dt + v(\rho',X_{\rho'},I_{\rho'})& \\
- \sum_{n\in\N}c(X_{\tilde\tau_n'},I_{(\tilde\tau_n')^-},I_{\tilde\tau_n'})1_{\{\tau'\leq\tilde\tau_n'<\rho'\}}\bigg)1_{\{(v^1-v^2)(\tau',X_{\tau'},I_{\tau'})<0\}}\bigg|\Fc_{\tau'}^s\bigg].&
\end{align*}
In conclusion, we find
\begin{align*}
&v(\tau',X_{\tau'},I_{\tau'}) \\
&= \ v^1(\tau',X_{\tau'},I_{\tau'})1_{\{(v^1-v^2)(\tau',X_{\tau'},I_{\tau'})\geq0\}} + v^2(\tau',X_{\tau'},I_{\tau'})1_{\{(v^1-v^2)(\tau',X_{\tau'},I_{\tau'})<0\}} \\
&\leq \ \E\bigg[\int_{\tau'}^{\rho'} f(X_t,I_t,u_t)dt + v(\rho',X_{\rho'},I_{\rho'}) - \sum_{n\in\N}c(X_{\tilde\tau_n'},I_{(\tilde\tau_n')^-},I_{\tilde\tau_n'})1_{\{\tau'\leq\tilde\tau_n'<\rho'\}}\bigg|\Fc_{\tau'}^s\bigg],
\end{align*}
which shows that $v\in\Vc^-$.

A similar argument allows to prove the stability with respect to infimum of $\Vc^+$ in (ii). In particular, fix $s\in[0,T]$, $\tau\in\Tc^s$, and $\alpha=(\tau_n,\iota_n)_{n\in\N}\in\Ac_s$. Let $\tilde u^1,\tilde u^2\in\Uc_{s,\tau}^E$ be the two elementary feedback strategies, for the nature, starting at $\tau$ and corresponding to $v^1$ and $v^2$. Let $B:=\{(v^1-v^2)(\tau(y),y(\tau(y)^+))\leq0\}\in\Bc_{\tau^+}^s$. Then, from Lemma \ref{L:Partition} we see that the map
\beqs
\tilde u(t,y) &:=& \tilde u^1(t,y) 1_{\{y\in B\}} + \tilde u^2(t,y) 1_{\{y\in B^c\}}
\enqs
is an elementary feedback strategy starting at $\tau$, which allows to prove that $v\in\Vc^+$.
\ep

\begin{Lemma}
\label{L:ApproximatingSequence}
Let Assumptions {\bf (H1)} and {\bf (H2)} hold.
\begin{itemize}
\item[\textup{(i)}] There exists a nondecreasing sequence $(v_n)_{n\in\N}\subset\Vc^-$ such that $v_n\nearrow v^-$.
\item[\textup{(ii)}] There exists a nonincreasing sequence $(v_n)_{n\in\N}\subset\Vc^+$ such that $v_n\searrow v^+$.
\end{itemize}
\end{Lemma}
\textbf{Proof.}
From Proposition 4.1 in \cite{BS12} we can find a sequence $(\tilde v_n)_{n\in\N}\subset\Vc^-$ satisfying $v^-=\sup_{n\in\N}\tilde v_n$. Set $v_n:=\tilde v_0\vee\cdots\vee\tilde v_n$, $n\in\N$. Then $v_n\nearrow v^-$ as $n\rightarrow\infty$, and from Lemma \ref{L:Stability} we see that $(v_n)_{n\in\N}\subset\Vc^-$. In a similar way we can prove statement (ii).
\ep

\vspace{3mm}

We are now in a position to prove Theorem \ref{T:StochPerron}. Firstly, we just state here, in the spirit of Lemma 2.4 in \cite{BS14}, the following technical result, which will be used several times in the proof of Theorem \ref{T:StochPerron}.

\begin{Lemma}
\label{L:Dini}
Let $\Cc\subset[0,T]\times\R^d$ be a compact set and consider a continuous function $F\colon\R^m\times\Cc\rightarrow\R$, which is nondecreasing in each of its first $m$ components. If there exists $\delta>0$ such that $\inf_{(t,x)\in\Cc}F(v^-(t,x,\cdot),t,x)>\delta$ $($resp. $\sup_{(t,x)\in\Cc}F(v^+(t,x,\cdot),t,x)<-\delta$$)$, then
\begin{align*}
&\inf_{(t,x)\in\Cc}F(v(t,x,\cdot),t,x) \ > \ \delta \\
\Big(\text{resp. }&\sup_{(t,x)\in\Cc}F(v(t,x,\cdot),t,x) \ < \ -\delta\Big)
\end{align*}
for some $v\in\Vc^-$ $($resp. $v\in\Vc^+$$)$.
\end{Lemma}
\textbf{Proof.}
Notice that, from the strict inequality $\inf_{(t,x)\in\Cc}F(v^-(t,x,\cdot),t,x)>\delta$ we can find $\eps>0$ such that $F(v^-(t,x,\cdot),t,x)>\delta+\eps$, for any $(t,x)\in\Cc$. Recall from Lemma \ref{L:ApproximatingSequence} that there exists a nondecreasing sequence $(v_n)_{n\in\N}\subset\Vc^-$ such that $v_n\nearrow v^-$. Let
\[
A_n \ := \ \big\{(t,x)\in\Cc\colon F(v_n(t,x,\cdot),t,x)\leq\delta+\eps/2\big\}.
\]
Notice that $A_n$ is closed, $A_{n+1}\subset A_n$, and $\cap_{n=0}^\infty A_n=\emptyset$. Since $A_n\subset\Cc$, using the compactness we see that there exists an $n_0$ such that $A_{n_0}=\emptyset$, namely $F(v_{n_0}(t,x,\cdot),t,x)>\delta+\eps$, for any $(t,x)\in\Cc$. In particular, $\inf_{(t,x)\in\Cc}F(v_{n_0}(t,x,\cdot),t,x)>\delta$. We then take $v:=v_{n_0}$. In a similar way we can prove the statement for $v^+$.
\ep

\vspace{3mm}

\noindent\textbf{Proof of Theorem \ref{T:StochPerron}.}\\
\textbf{Step I}. \emph{$v^-$ is a viscosity supersolution to the HJB equation \eqref{HJB}.}\\
\textbf{Step I(i)}. \emph{Interior viscosity supersolution property.} Let $(t_0,x_0)\in[0,T)\times\R^d$, $i\in\I_m$, and consider a test function $\varphi\in C^{1,2}([0,T]\times\R^d)$ such that $v^-(\cdot,\cdot,i)-\varphi(\cdot,\cdot)$ attains a strict global minimum equal to zero at $(t_0,x_0)$. Reasoning by contradiction, we assume that
\beqs
\min\Big\{-\dfrac{\partial \varphi}{\partial t}(t_0,x_0) - \inf_{u\in U}\big[\Lc^{i,u}\varphi(t_0,x_0) + f(x_0,i,u)\big], & & \\
v^-(t_0,x_0,i) - \max_{j\neq i}\big[v^-(t_0,x_0,j) - c(x_0,i,j)\big]\Big\} &<& 0.
\enqs
We distinguish two cases.

\noindent\textbf{Case a}. $-\frac{\partial \varphi}{\partial t}(t_0,x_0) - \inf_{u\in U}[\Lc^{i,u}\varphi(t_0,x_0) + f(x_0,i,u)]<0$. Then, there exists $\eps\in(0,T-t_0)$ such that
\[
-\frac{\partial \varphi}{\partial t}(t_0,x_0) - \inf_{u\in U}\big[\Lc^{i,u}\varphi(t_0,x_0) + f(x_0,i,u)\big] \ < \ -\eps.
\]
From the continuity of $b,\sigma,f$, together with the compactness of $U$, we see that we can choose a smaller $\eps\in(0,T-t_0)$ such that
\[
-\frac{\partial \varphi}{\partial t}(t,x) - \inf_{u\in U}\big[\Lc^{i,u}\varphi(t,x) + f(x,i,u)\big] \ < \ -\eps, \qquad \forall\,(t,x)\in B(t_0,x_0,\eps),
\]
where
\begin{equation}
\label{B(t_0,x_0)}
B(t_0,x_0,\eps) \ = \ \big\{(t,x)\in[0,T]\times\R^d\colon\max\{|t-t_0|,|x-x_0|\}<\eps\big\}.
\end{equation}
Since $v^-(\cdot,\cdot,i)-\varphi(\cdot,\cdot)$ is lsc and strictly positive on the compact set $\Cc:=\overline{B(t_0,x_0,\eps)}\backslash B(t_0,x_0,\eps/2)$, there exists $\delta>0$ such that $\inf_{(t,x)\in\Cc}(v^-(t,x,i)-\varphi(t,x))>\delta$. Denoting $F(p,t,x):=p-\varphi(t,x)$, it follows from Lemma \ref{L:Dini} that there exists $v\in\Vc^-$ such that $\varphi(t,x)+\delta<v(t,x,i)$ on $\Cc$. Now, define
\[
v^\delta(t,x,i) \ = \
\begin{cases}
(\varphi(t,x)+\delta)\vee v(t,x,i), \qquad &\text{on }\overline{B(t_0,x_0,\eps)}, \\
v(t,x,i), &\text{outside }\overline{B(t_0,x_0,\eps)}.
\end{cases}
\]
Moreover, $v^\delta(t,x,j)=v(t,x,j)$ for any $(t,x,j)\in[0,T]\times\R^d\times\I_m$, with $j\neq i$. Our aim is to prove that $v^\delta\in\Vc^-$, which would give a contradiction, since $v^\delta(t_0,x_0,i)>v^-(t_0,x_0,i)$. Clearly, $v^\delta$ satisfies the first item in Definition \ref{D:StochSubSol}, therefore it remains to prove the second item. To this end, fix $s\in[0,T]$ and $\tau,\rho\in\Tc^s$, with $\tau\leq\rho\leq T$. Let $\tilde\alpha^0=(\tilde\tau_n^0,\tilde\iota_n^0)_{n\in\N}$ be given by
\[
(\tilde\tau_n^0,\tilde\iota_n^0) \ = \ (T,i), \quad \forall\,n\in\N.
\]
Notice that $\tilde\alpha^0\in\Ac_{s,\tau^+}$. Introduce now the stopping rule $\rho_1\colon C([s,T];\R^d)\times\mathscr{L}([s,T];\I_m)\rightarrow[s,T]$, $\tau\leq\rho_1\leq T$,
\begin{equation}
\label{rho1}
\rho_1(y) \ = \ \inf\big\{t\in[\tau(y),T]\colon(t,y^X(t))\notin B(t_0,x_0,\eps/2)\big\}\wedge T.
\end{equation}
We denote by $\tilde\alpha^1=(\tilde\tau_n^1,\tilde\iota_n^1)_{n\in\N}\in\Ac_{s,(\rho_1\wedge\rho)^+}$ the feedback switching strategy in Definition \ref{D:StochSubSol}, corresponding to $s,\rho_1\wedge\rho,\rho$, for $v$. Then, we define $\tilde\alpha^2=\tilde\alpha^0\otimes_{\rho_1\wedge\rho}\tilde\alpha^1$, which belongs to $\Ac_{s,\tau^+}$  thanks to Proposition \ref{P:Concatenate}. Moreover, let $\tilde\alpha^3=(\tilde\tau_n^3,\tilde\iota_n^3)_{n\in\N}\in\Ac_{s,\tau^+}$ be the feedback switching strategy corresponding to $s,\tau,\rho$ for $v$. Then, we define $\tilde\alpha=(\tilde\tau_n,\tilde\iota_n)_{n\in\N}$ by (for any $y\in C([s,T];\R^d)\times\mathscr{L}([s,T];\I_m)$ we write $y=(y^X,y^I)$ with $y^X\in C([s,T];\R^d)$ and $y^I\in\mathscr{L}([s,T];\I_m)$)
\begin{align*}
(\tilde\tau_n(y),\tilde\iota_n(y)) \ &= \ (\tilde\tau_n^2(y),\tilde\iota_n^2(y)) 1_{\{(\tau(y),y^X(\tau(y)))\in B(t_0,x_0,\eps),\,(v-\varphi)(\tau(y),y(\tau(y)^+))<\delta,\,y^I(\tau(y)^+)=i\}} \\
&\quad \ + (\tilde\tau_n^3(y),\tilde\iota_n^3(y)) 1_{\{(\tau(y),y^X(\tau(y)))\in B(t_0,x_0,\eps),\,(v-\varphi)(\tau(y),y(\tau(y)^+))<\delta,\,y^I(\tau(y)^+)=i\}^c}.
\end{align*}
From Lemma \ref{L:Partition} it follows that $\tilde\alpha\in\Ac_{s,\tau^+}$. Moreover, the feedback switching strategy $\tilde\alpha$ satisfies the condition in the second item of Definition \ref{D:StochSubSol} for $v^\delta$. To see this, fix $\alpha=(\tau_n,\iota_n)_{n\in\N}\in\Ac_s$, $u\in\Uc_s$, and $(x,i)\in\R^d\times\I_m$. We adopt the shorthands:
\beqs
(X,I) &=& (X^{s,x,i;\alpha\otimes_\tau\tilde\alpha,u},I^{s,x,i;\alpha\otimes_\tau\tilde\alpha,u}), \\
(X^1,I^1) &=& (X^{s,x,i;\alpha\otimes_\tau\tilde\alpha^2,u},I^{s,x,i;\alpha\otimes_\tau\tilde\alpha^2,u}), \\
(X^2,I^2) &=& (X^{s,x,i;\alpha\otimes_\tau\tilde\alpha^3,u},I^{s,x,i;\alpha\otimes_\tau\tilde\alpha^3,u}).
\enqs
We also denote $\tau'=\tau(X_\cdot,I_{\cdot^-})$, $\rho_1'=\rho_1(X_\cdot,I_{\cdot^-})$, and $\rho'=\rho(X_\cdot,I_{\cdot^-})$. Notice that
\begin{align*}
(X,I) \ &= \ (X^1,I^1)1_{\{(\tau',X_{\tau'})\in B(t_0,x_0,\eps),\,(v-\varphi)(\tau',X_{\tau'},I_{\tau'})<\delta,\,I_{\tau'}=i\}} \\
&\quad \ + (X^2,I^2)1_{\{(\tau',X_{\tau'})\in B(t_0,x_0,\eps),\,(v-\varphi)(\tau',X_{\tau'},I_{\tau'})<\delta,\,I_{\tau'}=i\}^c}.
\end{align*}
In particular, it is useful to decompose $v^\delta(\tau',X_{\tau'},I_{\tau'})$ as follows
\begin{align}
\label{v^delta_Decomposition}
v^\delta(\tau',X_{\tau'},I_{\tau'}) \ &= \ \big(\varphi(\tau',X_{\tau'}^1) + \delta\big) 1_{\{(\tau',X_{\tau'})\in B(t_0,x_0,\eps),\,(v-\varphi)(\tau',X_{\tau'},I_{\tau'})<\delta,\,I_{\tau'}=i\}} \\
&\quad \ + v(\tau',X_{\tau'}^2,I_{\tau'}^2) 1_{\{(\tau',X_{\tau'})\in B(t_0,x_0,\eps),\,(v-\varphi)(\tau',X_{\tau'},I_{\tau'})<\delta,\,I_{\tau'}=i\}^c}. \notag
\end{align}
We now consider the two terms on the right-hand side of \eqref{v^delta_Decomposition} individually. Regarding the first term, we apply It\^o's formula to $\varphi$ between $\tau'$ and $\rho_1'\wedge\rho'$, observing that $I_t^1=i$ for any $t\in[\tau',\rho_1'\wedge\rho']$; afterwards, we use the property in the second item of Definition \ref{D:StochSubSol} for $v$ with corresponding feedback switching strategy $\tilde\alpha^1$. Finally, concerning the other term in \eqref{v^delta_Decomposition}, the result follows from the properties of $v$ and the definition of $\tilde\alpha^3$.

\noindent\textbf{Case b}. $v^-(t_0,x_0,i)<\max_{j\neq i}[v^-(t_0,x_0,j) - c(x_0,i,j)]$ and $-\frac{\partial \varphi}{\partial t}(t_0,x_0) - \inf_{u\in U}[\Lc^{i,u}\varphi(t_0,x_0) + f(x_0,i,u)]\geq0$. Since $v^-$ is lsc and $c$ is continuous, there exists $\eps\in(0,T-t_0)$ such that
\[
v^-(t_0,x_0,i) + \eps \ < \ \inf_{(t,x)\in\overline{B(t_0,x_0,\eps)}}\max_{j\neq i}[v^-(t,x,j) - c(x,i,j)].
\]
Set $F(p,t,x)=\max_{j\neq i}[p_j - c(x,i,j)]$, for any $(p,t,x)\in\R^m\times\overline{B(t_0,x_0,\eps)}$. Then, from Lemma \ref{L:Dini} it follows that there exists $v\in\Vc^-$ such that $F(v(t,x,\cdot),t,x)>v^-(t_0,x_0,i)+\eps\geq v(t_0,x_0,i)+\eps$, for any $(t,x)\in\overline{B(t_0,x_0,\eps)}$. We also suppose that the function $v$ given by Lemma \ref{L:Dini} satisfies $v^-(t_0,x_0,i)-v(t_0,x_0,i)<\eps/2$. Since $v$ is continuous on $\overline{B(t_0,x_0,\eps)}$, we can find $\delta>0$ such that
\begin{equation}
\label{sup_v+eps<max_v-c}
\sup_{(t',x')\in\overline{B(t_0,x_0,\delta)}}v(t',x',i) + \eps \ < \ \inf_{(t,x)\in\overline{B(t_0,x_0,\eps)}}\max_{j\neq i}\big[v(t,x,j) - c(x,i,j)\big].
\end{equation}
Let $M>0$ be an upper bound for the continuous function $|f(x,i,u)|$ on the compact set $\overline{B(t_0,x_0,\eps)}\times\I_m\times U$. We suppose that $\delta\leq\eps/(4M)$. Now, define (we adopt the notation $\|(t,x)\|=\max\{|t|,|x|\}$)
\[
v^\delta(t,x,i) \ = \
\begin{cases}
v(t,x,i) + \frac{\eps}{2\delta}(\delta - \|(t-t_0,x-x_0)\|), \qquad &\text{on }\overline{B(t_0,x_0,\delta)}, \\
v(t,x,i), &\text{outside }\overline{B(t_0,x_0,\delta)}.
\end{cases}
\]
Moreover, $v^\delta(t,x,j)=v(t,x,j)$ for any $(t,x,j)\in[0,T]\times\R^d\times\I_m$, with $j\neq i$. As $v^\delta(t_0,x_0,i)>v^-(t_0,x_0,i)$, we get a contradiction if we prove that $v^\delta\in\Vc^-$. In order to do so, fix $s\in[0,T]$ and $\tau,\rho\in\Tc^s$, with $\tau\leq\rho\leq T$. We have to determine $\tilde\alpha=(\tilde\tau_n,\tilde\iota_n)_{n\in\N}\in\Ac_{s,\tau^+}$ which works for $v^\delta$. To this end, define $\rho_1\in\Tc^s$ as follows
\[
\rho_1(y) \ = \ \inf\big\{t\in[\tau(y),T]\colon(t,y^X(t))\notin B(t_0,x_0,\delta)\big\}\wedge T.
\]
Let $\tilde\alpha^0=(\tilde\tau_n^0,\tilde\iota_n^0)_{n\in\N}$ be given by: $(\tilde\tau_n^0,\tilde\iota_n^0)=(T,i)$ for any $n\geq1$, and
\begin{align*}
\tilde\tau_0^0(y) \ &= \ \rho_1(y) 1_{\{(\tau(y),y^X(\tau(y)))\in B(t_0,x_0,\delta)\}} + T\,1_{\{(\tau(y),y^X(\tau(y)))\notin B(t_0,x_0,\delta)\}}, \\
\tilde\iota_0^0(y) \ &= \ \min\big\{j\neq i\colon v(\tilde\tau_0^0(y),y^X(\tilde\tau_0^0(y)),j) - c(y^X(\tilde\tau_0^0(y)),i,j)=m(y)\big\},
\end{align*}
where $m\colon C([s,T];\R^d)\times\mathscr{L}([s,T];\I_m)\rightarrow\R$ is defined as
\[
m(y) \ = \ \max_{j\neq i}\big[v(\tilde\tau_0^0(y),y^X(\tilde\tau_0^0(y)),j) - c(y^X(\tilde\tau_0^0(y)),i,j)\big].
\]
Notice that $m$ is $\Bc_{\tilde\tau_0^0}^s$-measurable, so that $\tilde\iota_0^0$ is $\Bc_{\tilde\tau_0^0}^s$-measurable. Moreover, $\tau<\tilde\tau_0^0$ on the set $\{\tau<T\}$. In particular, $\tilde\alpha^0\in\Ac_{\tau^+}^s$. Now, consider the feedback switching strategy $\tilde\alpha^1=(\tilde\tau_n^1,\tilde\iota_n^1)_{n\in\N}\in\Ac_{s,(\tilde\tau_0^0\wedge\rho)^+}$ in Definition \ref{D:StochSubSol}, corresponding to $s,\tilde\tau_0^0\wedge\rho,\rho$, for $v$. We define $\tilde\alpha^2=\tilde\alpha^0\otimes_{\tilde\tau_0^0\wedge\rho}\tilde\alpha^1$, which belongs to $\Ac_{s,\tau^+}$ thanks to Proposition \ref{P:Concatenate}. Consider also the feedback switching strategy $\tilde\alpha^3=(\tilde\tau_n^3,\tilde\iota_n^3)_{n\in\N}\in\Ac_{s,\tau^+}$, corresponding to $s,\tau,\rho$, for $v$. Then, let $\tilde\alpha=(\tilde\tau_n,\tilde\iota_n)_{n\in\N}$ be given by
\begin{align*}
(\tilde\tau_n(y),\tilde\iota_n(y)) \ &= \ (\tilde\tau_n^2(y),\tilde\iota_n^2(y)) 1_{\{(\tau(y),y^X(\tau(y)))\in B(t_0,x_0,\delta),\,y^I(\tau(y)^+)=i\}} \\
&\quad \ + (\tilde\tau_n^3(y),\tilde\iota_n^3(y)) 1_{\{(\tau(y),y^X(\tau(y)))\in B(t_0,x_0,\delta),\,y^I(\tau(y)^+)=i\}^c}.
\end{align*}
From Lemma \ref{L:Partition} it follows that $\tilde\alpha\in\Ac_{s,\tau^+}$. Moreover, $\tilde\alpha$ is the feedback switching strategy which satisfies the condition in the second item of Definition \ref{D:StochSubSol} for $v^\delta$. To see this, fix $\alpha=(\tau_n,\iota_n)_{n\in\N}\in\Ac_s$, $u\in\Uc_s$, and $(x,i)\in\R^d\times\I_m$. We adopt the shorthands introduced in Case $a$. Consider the event $A:=\{(\tau',X_{\tau'})\in B(t_0,x_0,\delta),\,I_{\tau'}=i\}$. On $A^c$ the result follows from the properties of $v$ and the definition of $\tilde\alpha^3$. On the other hand, on $A$ we have
\begin{align*}
v^\delta(\tau',X_{\tau'},I_{\tau'}) 1_A \ = \ v^\delta(\tau',X_{\tau'}^1,i) 1_A \ &= \ \Big[v(\tau',X_{\tau'}^1,i) + \frac{\eps}{2\delta}\big(\delta - \|(\tau'-t_0,X_{\tau'}^1-x_0)\|\big)\Big] 1_A \\
&\leq \ \Big[v(\tau',X_{\tau'}^1,i) + \frac{\eps}{2}\Big] 1_A.
\end{align*}
Using \eqref{sup_v+eps<max_v-c} and taking the conditional expectation with respect to $\Fc_{\tau'}^s$, we obtain (denoting $\tilde\tau_0^{0,'}=\tilde\tau_0^0(X_\cdot,I_{\cdot^-})$)
\[
v^\delta(\tau',X_{\tau'},I_{\tau'}) 1_A \ \leq \ \E\Big[v\big(\tilde\tau_0^{0,'}\wedge\rho',X_{\tilde\tau_0^{0,'}\wedge\rho'}^1,I_{\tilde\tau_0^{0,'}\wedge\rho'}^1\big) - c\big(X_{\tilde\tau_0^{0,'}\wedge\rho'}^1,i,I_{\tilde\tau_0^{0,'}\wedge\rho'}^1\big) - \frac{\eps}{2}\Big|\Fc_{\tau'}^s\Big] 1_A.
\]
Observe that $\tilde\tau_0^{0,'}\leq\rho'$ on $A$. Therefore, the above inequality can be written as
\[
v^\delta(\tau',X_{\tau'},I_{\tau'}) 1_A \ \leq \ \E\Big[v\big(\tilde\tau_0^{0,'},X_{\tilde\tau_0^{0,'}}^1,I_{\tilde\tau_0^{0,'}}^1\big) - c\big(X_{\tilde\tau_0^{0,'}}^1,i,I_{\tilde\tau_0^{0,'}}^1\big) - \frac{\eps}{2}\Big|\Fc_{\tau'}^s\Big] 1_A.
\]
Adding and subtracting $\int_{\tau'}^{\tilde\tau_0^{0,'}}f(X_t^1,I_t^1,u_t)dt$, noting that $(\tilde\tau_0^{0,'}-\tau')1_A\leq2\delta$ and $2\delta M-\eps/2\leq0$, we find
\[
v^\delta(\tau',X_{\tau'},I_{\tau'}) 1_A \, \leq \, \E\bigg[\int_{\tau'}^{\tilde\tau_0^{0,'}}f(X_t^1,I_t^1,u_t)dt + v\big(\tilde\tau_0^{0,'},X_{\tilde\tau_0^{0,'}}^1,I_{\tilde\tau_0^{0,'}}^1\big) - c\big(X_{\tilde\tau_0^{0,'}}^1,i,I_{\tilde\tau_0^{0,'}}^1\big)\bigg|\Fc_{\tau'}^s\bigg] 1_A.
\]
Finally, using that $v$ satisfies the second item of Definition \ref{D:StochSubSol}, with corresponding feedback switching strategy $\tilde\alpha^1$, and from the inequality $v\leq v^\delta$, we deduce that $v^\delta\in\Vc^-$.

\noindent\textbf{Step I(ii)}. \emph{Terminal condition.} Reasoning by contradiction, we assume that there exist $x_0\in\R^d$ and $i\in\I_m$ such that
\[
v^-(T,x_0,i) \ < \ g(x_0,i).
\]
Since $g$ is continuous, there exists $\eps>0$ such that $v^-(T,x_0,i)\leq g(x,i)-\eps$ whenever $|x-x_0|\leq\eps$. Consider the compact set
\[
\Cc \ := \ \big(\overline{B(T,x_0,\eps)}\backslash B(T,x_0,\eps/2)\big)\cap\big([0,T]\times\R^d\big),
\]
where $B(T,x_0,\eps)=\{(t,x)\in[0,T]\times\R^d\colon\max\{|t-t_0|,|x-x_0|\}<\eps\}$. Since $v^-$ is lsc, it is bounded from below on $\Cc$. Therefore, we can find $\eta>0$ small enough (possibly depending on $\eps$) such that
\[
v^-(T,x_0,i) - \frac{\eps^2}{4\eta} \ < \ -\eps + \inf_{(t,x)\in\Cc}v^-(t,x,i).
\]
From Lemma \ref{L:Dini} with $F(p,t,x)=p$ for any $(p,t,x)\in\R\times\Cc$, we can find $v\in\Vc^-$ such that
\begin{equation}
\label{v^-<inf_v}
v^-(T,x_0,i) - \frac{\eps^2}{4\eta} \ < \ -\eps + \inf_{(t,x)\in\Cc}v(t,x,i).
\end{equation}
For $k>0$ define
\[
\varphi^{\eta,\eps,k}(t,x) \ = \ v^-(T,x_0,i) - \frac{|x-x_0|^2}{\eta} - k(T-t).
\]
Since $b,\sigma,f$ are continuous, we can choose $k$ large enough such that
\[
-\frac{\partial \varphi^{\eta,\eps,k}}{\partial t}(t,x) - \inf_{u\in U}\big[\Lc^{i,u}\varphi^{\eta,\eps,k}(t,x) + f(x,i,u)\big] \ < \ 0, \qquad \forall\,(t,x)\in\overline{B(T,x_0,\eps)}.
\]
From \eqref{v^-<inf_v} it follows that $\varphi^{\eta,\eps,k}(t,x)<-\eps+v(t,x,i)$ on $\Cc$. Moreover
\[
\varphi^{\eta,\eps,k}(T,x) \ \leq \ v^-(T,x_0,i) \ \leq \ g(x,i) - \eps, \qquad \text{whenever }|x-x_0|\leq\eps.
\]
Now, for $\delta\in(0,\eps)$ define
\[
v^\delta(t,x,i) \ = \
\begin{cases}
(\varphi^{\eta,\eps,k}(t,x)+\delta)\vee v(t,x,i), \qquad &\text{on }\overline{B(t_0,x_0,\eps)}, \\
v(t,x,i), &\text{outside }\overline{B(t_0,x_0,\eps)}.
\end{cases}
\]
Moreover, $v^\delta(t,x,j)=v(t,x,j)$ for any $(t,x,j)\in[0,T]\times\R^d\times\I_m$, with $j\neq i$. As $v^\delta(T,x_0,i)>v^-(T,x_0,i)$, we get a contradiction if we are able to prove that $v^\delta\in\Vc^-$. In particular, for any $s\in[0,T]$ and $\tau,\rho\in\Tc^s$ with $\tau\leq\rho\leq T$, we have to find $\tilde\alpha=(\tilde\tau_n,\tilde\iota_n)_{n\in\N}\in\Ac_{s,\tau^+}$ which works for $v^\delta$. Consider the feedback switching strategy $\tilde \alpha$ defined in Step I(i), Case $a$, with $\rho_1$ the exit time from $B(T,x_0,\eps/2)$. Then, proceeding as in Case $a$ of Step I(i), we can prove that $\tilde\alpha$ satisfies the condition in the second item of Definition \ref{D:StochSubSol} for $v^\delta$.

\noindent\textbf{Step II}. \emph{$v^+$ is a viscosity subsolution to the HJB equation \eqref{HJB}.}\\
\textbf{Step II(i)}. \emph{Interior viscosity subsolution property.} Let $(t_0,x_0)\in[0,T)\times\R^d$, $i\in\I_m$, and consider a test function $\varphi\in C^{1,2}([0,T]\times\R^d)$ such that $v^+(\cdot,\cdot,i)-\varphi(\cdot,\cdot)$ attains a strict global maximum equal to zero at $(t_0,x_0)$. Reasoning by contradiction, we assume that
\beqs
\min\Big\{-\dfrac{\partial \varphi}{\partial t}(t_0,x_0) - \inf_{u\in U}\big[\Lc^{i,u}\varphi(t_0,x_0) + f(x_0,i,u)\big], & & \\
v^+(t_0,x_0,i) - \max_{j\neq i}\big[v^+(t_0,x_0,j) - c(x_0,i,j)\big]\Big\} &>& 0.
\enqs
Then, there exists $\eps>0$ and $\underline u\in U$ such that
\[
-\frac{\partial \varphi}{\partial t}(t_0,x_0) - \Lc^{i,\underline u}\varphi(t_0,x_0) - f(x_0,i,\underline u) \ > \ \eps.
\]
From the continuity of $b,\sigma,f$, it follows that we can find a smaller $\eps>0$ such that
\[
-\frac{\partial \varphi}{\partial t}(t,x) - \Lc^{i,\underline u}\varphi(t,x) - f(x,i,\underline u) \ > \ \eps, \qquad \forall\,(t,x)\in B(t_0,x_0,\eps),
\]
where $B(t_0,x_0,\eps)$ is given by \eqref{B(t_0,x_0)}. As $v^+(\cdot,\cdot,i)-\varphi(\cdot,\cdot)$ is usc and strictly negative on the compact set $\Cc:=\overline{B(t_0,x_0,\eps)}\backslash B(t_0,x_0,\eps/2)$, we see that there exists $\delta>0$ such that $\sup_{(t,x)\in\Cc}(v^+(t,x,i)-\varphi(t,x))<-\delta$. Denoting $F(p,t,x):=p-\varphi(t,x)$, it follows from Lemma \ref{L:Dini} that there exists $v\in\Vc^+$ such that $\varphi(t,x)-\delta>v(t,x,i)$ on $\Cc$. Now, define
\[
v^\delta(t,x,i) \ = \
\begin{cases}
(\varphi(t,x)-\delta)\wedge v(t,x,i), \qquad &\text{on }\overline{B(t_0,x_0,\eps)}, \\
v(t,x,i), &\text{outside }\overline{B(t_0,x_0,\eps)}.
\end{cases}
\]
Moreover, $v^\delta(t,x,j)=v(t,x,j)$ for any $(t,x,j)\in[0,T]\times\R^d\times\I_m$, with $j\neq i$. As $v^\delta(t_0,x_0,i)<v^+(t_0,x_0,i)$, we find a contradiction if we are able to prove that $v^\delta\in\Vc^+$. To this end, fix $s\in[0,T]$, $\tau\in\Tc^s$, and $\alpha=(\tau_n,\iota_n)_{n\in\N}\in\Ac_s$. We have to construct an elementary feedback strategy  $\tilde u\in\Uc_{s,\tau}^E$ which works for $v^\delta$. Consider the stopping rule $\rho_1\in\Tc^s$ given by \eqref{rho1}, and let $\tilde u^1\in\Uc_{s,\rho_1}^E$ be the elementary feedback strategy for $v$, corresponding to $s,\rho_1,\alpha$. Then, we define $\tilde u^2=\underline u\otimes_{\rho_1}\tilde u^1$, which belongs to $\Uc_{s,\tau}^E$ thanks to Proposition \ref{P:Concatenate}. Now, let $\tilde u^3\in\Uc_{s,\tau}^E$ be the elementary feedback strategy for $v$, corresponding to $s,\tau,\alpha$. Then, we define
\begin{align*}
\tilde u(t,y) \ &= \ \tilde u^2(t,y) 1_{\{(\tau(y),y^X(\tau(y)))\in B(t_0,x_0,\eps),\,(v-\varphi)(\tau(y),y(\tau(y)^+))>-\delta,\,y^I(\tau(y)^+)=i\}} \\
&\quad \ + \tilde u^3(t,y) 1_{\{(\tau(y),y^X(\tau(y)))\in B(t_0,x_0,\eps),(v-\varphi)(\tau(y),y(\tau(y)^+))>-\delta\,,\,y^I(\tau(y)^+)=i\}^c}.
\end{align*}
From Lemma \ref{L:Partition} we see that $\tilde u\in\Uc_{s,\tau}^E$. Moreover, $\tilde u$ is the elementary feedback strategy for the second item of Definition \ref{D:StochSuperSol} for $v^\delta$. Indeed, fix $u\in\Uc_s^E$, $(x,i)\in\R^d\times\I_m$, and $\rho\in\Tc^s$, with $\tau\leq\rho\leq T$. We adopt the shorthands:
\beqs
(X,I) &=& (X^{s,x,i;\alpha,u\otimes_\tau\tilde u},I^{s,x,i;\alpha,u\otimes_\tau\tilde u}), \\
(X^1,I^1) &=& (X^{s,x,i;\alpha,u\otimes_\tau\tilde u^2},I^{s,x,i;\alpha,u\otimes_\tau\tilde u^2}), \\
(X^2,I^2) &=& (X^{s,x,i;\alpha,u\otimes_\tau\tilde u^3},I^{s,x,i;\alpha,u\otimes_\tau\tilde u^3}).
\enqs
We also denote $\tau'=\tau(X_\cdot,I_{\cdot^-})$, $\rho_1'=\rho_1(X_\cdot,I_{\cdot^-})$, and $\rho'=\rho(X_\cdot,I_{\cdot^-})$. Notice that
\begin{align*}
(X,I) \ &= \ (X^1,I^1) 1_{\{(\tau',X_{\tau'})\in B(t_0,x_0,\eps),\,(v-\varphi)(\tau',X_{\tau'},I_{\tau'})>-\delta,\,I_{\tau'}=i\}} \\
&\quad \ + (X^2,I^2) 1_{\{(\tau',X_{\tau'})\in B(t_0,x_0,\eps),\,(v-\varphi)(\tau',X_{\tau'},I_{\tau'})>-\delta,\,I_{\tau'}=i\}^c}.
\end{align*}
Moreover, write $v^\delta(\tau',X_{\tau'},I_{\tau'})$ as follows
\begin{align*}
v^\delta(\tau',X_{\tau'},I_{\tau'}) \ &= \ \big(\varphi(\tau',X_{\tau'}^1) - \delta\big) 1_{\{(\tau',X_{\tau'})\in B(t_0,x_0,\eps),\,(v-\varphi)(\tau',X_{\tau'},I_{\tau'})>-\delta,\,I_{\tau'}=i\}} \\
&\quad \ + v(\tau',X_{\tau'}^2,I_{\tau'}^2) 1_{\{(\tau',X_{\tau'})\in B(t_0,x_0,\eps),\,(v-\varphi)(\tau',X_{\tau'},I_{\tau'})>-\delta,\,I_{\tau'}=i\}^c}.
\end{align*}
Then, applying It\^o's formula to $\varphi$ and using the properties of $v$, we see that $v^\delta\in\Vc^+$.

\noindent\textbf{Step II(ii)}. \emph{Terminal condition.} Reasoning by contradiction, we assume that there exist $x_0\in\R^d$ and $i\in\I_m$ such that
\[
v^+(T,x_0,i) \ > \ g(x_0,i).
\]
Since $g$ is continuous, there exists $\eps>0$ such that $v^+(T,x_0,i)\geq g(x,i)+\eps$ whenever $|x-x_0|\leq\eps$. Consider the compact set
\[
\Cc \ := \ \big(\overline{B(T,x_0,\eps)}\backslash B(T,x_0,\eps/2)\big)\cap\big([0,T]\times\R^d\big).
\]
As $v^+$ is usc, it is bounded from above on $\Cc$. Therefore, we can find $\eta>0$ small enough (possibly depending on $\eps$) such that
\[
v^+(T,x_0,i) + \frac{\eps^2}{4\eta} \ > \ \eps + \sup_{(t,x)\in\Cc}v^+(t,x,i).
\]
From Lemma \ref{L:Dini} with $F(p,t,x)=p$ for any $(p,t,x)\in\R\times\Cc$, we can find $v\in\Vc^+$ such that
\begin{equation}
\label{v^+>sup_v}
v^+(T,x_0,i) + \frac{\eps^2}{4\eta} \ > \ \eps + \sup_{(t,x)\in\Cc}v(t,x,i).
\end{equation}
For $k>0$ define
\[
\varphi^{\eta,\eps,k}(t,x) \ = \ v^+(T,x_0,i) + \frac{|x-x_0|^2}{\eta} + k(T-t).
\]
Since $b,\sigma,f$ are continuous, we can choose $k$ large enough and $\underline u\in U$ such that
\[
-\frac{\partial \varphi^{\eta,\eps,k}}{\partial t}(t,x) - \Lc^{i,\underline u}\varphi^{\eta,\eps,k}(t,x) - f(x,i,\underline u) \ > \ 0, \qquad \forall\,(t,x)\in\overline{B(T,x_0,\eps)}.
\]
From \eqref{v^+>sup_v} it follows that $\varphi^{\eta,\eps,k}(t,x)>\eps+v(t,x,i)$ on $\Cc$. Moreover
\[
\varphi^{\eta,\eps,k}(T,x) \ \geq \ v^+(T,x_0,i) \ \geq \ g(x,i) + \eps, \qquad \text{whenever }|x-x_0|\leq\eps.
\]
Now, for $\delta\in(0,\eps)$ define
\[
v^\delta(t,x,i) \ = \
\begin{cases}
(\varphi^{\eta,\eps,k}(t,x)-\delta)\wedge v(t,x,i), \qquad &\text{on }\overline{B(t_0,x_0,\eps)}, \\
v(t,x,i), &\text{outside }\overline{B(t_0,x_0,\eps)}.
\end{cases}
\]
Moreover, $v^\delta(t,x,j)=v(t,x,j)$ for any $(t,x,j)\in[0,T]\times\R^d\times\I_m$, with $j\neq i$. As $v^\delta(T,x_0,i)<v^+(T,x_0,i)$, we get a contradiction if we prove that $v^\delta\in\Vc^+$. In particular, for any $s\in[0,T]$, $\tau\in\Tc^s$, and $\alpha=(\tau_n,\iota_n)_{n\in\N}\in\Ac_s$, we have to find $\tilde u\in\Uc_{s,\tau}^E$  for the second item of Definition \ref{D:StochSuperSol} for $v^\delta$. Let $\tilde u\in\Uc_{s,\tau}^E$ be the elementary feedback strategy defined in Step II(i), with $\rho_1$ the exit time from $B(T,x_0,\eps/2)$. Then, we can prove, as in Step II(i), that $\tilde u$ satisfies the condition in the second item of Definition \ref{D:StochSuperSol} for $v^\delta$.
\ep

\section{Dynamic programming and viscosity properties of $V$}

\setcounter{equation}{0}
\setcounter{Theorem}{0} \setcounter{Proposition}{0}
\setcounter{Corollary}{0} \setcounter{Lemma}{0}
\setcounter{Definition}{0} \setcounter{Remark}{0}

In the present section, by means of the comparison principle for equation \eqref{HJB}, we prove that $V$ satisfies the dynamic programming principle and is a viscosity solution to equation \eqref{HJB}, which therefore turns out to be the dynamic programming equation of the robust switching control problem.

\subsection{Comparison principle and viscosity characterization}

We need to make an additional assumption on the switching costs in order to get comparison principle. 

\vspace{1mm}

{\bf (H3)}  
\begin{itemize}
\item[] The switching cost function $c$ satisfies the {\it no free loop property}:  for any sequence of indices $i_1,\ldots,i_k\in\I_m$, with $k\in\N\backslash\{0,1,2\}$,  $i_1=i_k$, and $\text{card}\{i_1,\ldots,i_k\}=k-1$, we have
\beqs
c(x,i_1,i_2) + c(x,i_2,i_3) + \cdots + c(x,i_{k-1},i_k) + c(x,i_k,i_1) \ > \ 0, \qquad \forall\,x\in\R^d.
\enqs
We also assume that $c(x,i,i)=0$, for any $x\in\R^d$ and $i\in\I_m$.
\end{itemize}

\begin{Theorem}[Comparison principle]
\label{CompThm}
Let Assumptions {\bf (H1)}, {\bf (H2)} and {\bf (H3)}  hold and consider a viscosity subsolution $\check v$ (resp. supersolution $\hat v$) to equation \eqref{HJB}. Suppose that
\beqs
\sup_{(t,x,i)\in[0,T]\times\R^d\times\I_m}\frac{|\check v(t,x,i)| + |\hat v(t,x,i)|}{1 + |x|^q} &<& \infty,
\enqs
for some $q\geq1$. Then, we have $\check v(t,x,i)\leq \hat v(t,x,i)$ for any $(t,x,i)\in[0,T]\times\R^d\times\I_m$.
\end{Theorem}

\begin{Remark}
{\rm
The proof can be done along the lines of Proposition 3.1 in \cite{ham_morl13}, apart from minor changes due to the presence of the infimum over $U$ in \eqref{HJB}, which are dealt with by the uniform Lipschitz condition in {\bf (H1)}(ii). 
More precisely, it is proved, as usual, proceeding by contradiction and then using the doubling variable technique. We simply notice here that equation \eqref{HJB} requires a particular step. Indeed, along the sequence of maximum points $(t_n,x_n)_n$ coming through the doubling of variables, we require
\beq
\label{CompIneq}
\check v(t_n,x_n,i) &>& \max_{j\neq i}\big[\check v(t_n,x_n,j) - c(x_n,i,j)\big],
\enq
so that, from the viscosity subsolution property of $\check v$, we can derive an inequality for the PDE part of equation \eqref{HJB} (concerning $\hat v$, the viscosity supersolution property implies already the nonnegativity of both terms in \eqref{HJB}). Condition \eqref{CompIneq} is obtained from a ``no-loop'' argument presented in Theorem 3.1 of \cite{ishii_koike91} (see also Lemma A.2 in \cite{barles_jak05} and Proposition 3.1 in \cite{ham_morl13}), which is based on the no free loop property in {\bf (H3)}. 
\ep
}
\end{Remark}

\begin{Corollary}
Under Assumptions {\bf (H1)}, {\bf (H2)}, and {\bf (H3)}, we have $v^-=V=\overline V=v^+$. In particular, $V$ (as $v^-,\overline V,v^+$) is continuous. Moreover, $V$ is the unique viscosity solution to equation \eqref{HJB} satisfying a polynomial growth condition. Furthermore, $V$ satisfies the dynamic programming principle: for any $(s,x,i)\in[0,T]\times\R^d\times\I_m$ and $\rho\in\Tc^s$,
\beqs
V(s,x,i) &=& \sup_{\alpha\in\Ac_{s^+}}\inf_{u\in\Uc_s} \E\bigg[\int_s^{\rho'} f(X_t,I_t,u_t)dt + V(\rho',X_{\rho'},I_{\rho'}) \\
& & \hspace{2.5cm} - \; \sum_{n\in\N}c(X_{\tau_n'},I_{(\tau_n')^-},I_{\tau_n'})1_{\{s\leq\tau_n'<\rho'\}}\bigg] \\
&=& \sup_{\alpha\in\Ac_{s^+}}\inf_{u\in\Uc_s^E} \E\bigg[\int_s^{\rho'} f(X_t,I_t,u_t')dt + V(\rho',X_{\rho'},I_{\rho'}) \\
& & \hspace{2.5cm} - \; \sum_{n\in\N}c(X_{\tau_n'},I_{(\tau_n')^-},I_{\tau_n'})1_{\{s\leq\tau_n'<\rho'\}}\bigg], \notag
\enqs
with the shorthands $X=X^{s,x,i;\alpha,u}$, $I=I^{s,x,i;\alpha,u}$, $\rho'=\rho(X_\cdot,I_{\cdot^-})$, $\tau_n'=\tau_n(X_\cdot,I_{\cdot^-})$, and $u_t'=u(t,X_\cdot,I_{\cdot^-})$.
\end{Corollary}
\textbf{Proof.}
The equality $v^-=V=\overline V=v^+$ follows from the comparison Theorem \ref{CompThm}. Since $v^-$ is lsc and $v^+$ is usc, we see that $V$ is continuous. Moreover, from Remark \ref{R:V_polynomial_growth}, Theorem \ref{T:StochPerron}, and Theorem \ref{CompThm} it follows that $V$ is the unique viscosity solution to equation \eqref{HJB} satisfying a polynomial growth condition. Finally, let us prove the dynamic programming principle for $V$. We begin noting that $v^-$ and $v^+$ satisfy, respectively, the sub- and super-dynamic programming principles: for any $(s,x,i)\in[0,T]\times\R^d\times\I_m$ and $\rho\in\Tc^s$,
\beq
\label{DPP^-_v^-}
v^-(s,x,i) &\leq & \sup_{\alpha\in\Ac_{s^+}}\inf_{u\in\Uc_s} \E\bigg[\int_s^{\rho'} f(X_t,I_t,u_t)dt + v^-(\rho',X_{\rho'},I_{\rho'}) \\
& & \hspace{2.5cm} - \; \sum_{n\in\N}c(X_{\tau_n'},I_{(\tau_n')^-},I_{\tau_n'})1_{\{s\leq\tau_n'<\rho'\}}\bigg] \notag
\enq
and
\beq
\label{DPP^+_v^+}
v^+(s,x,i) &\geq & \sup_{\alpha\in\Ac_{s^+}}\inf_{u\in\Uc_s^E} \E\bigg[\int_s^{\rho'} f(X_t,I_t,u_t')dt + v^+(\rho',X_{\rho'},I_{\rho'}) \\
& & \hspace{2.5cm} - \; \sum_{n\in\N}c(X_{\tau_n'},I_{(\tau_n')^-},I_{\tau_n'})1_{\{s\leq\tau_n'<\rho'\}}\bigg], \notag
\enq
with the shorthands $X=X^{s,x,i;\alpha,u}$, $I=I^{s,x,i;\alpha,u}$, $\rho'=\rho(X_\cdot,I_{\cdot^-})$, $\tau_n'=\tau_n(X_\cdot,I_{\cdot^-})$, and $u_t'=u(t,X_\cdot,I_{\cdot^-})$. As a matter of fact, let $(v_n)_{n\in\N}\subset\Vc^-$ be the sequence in Lemma \ref{L:ApproximatingSequence}(i). From Lemma \ref{L:V^-nonempty_HalfDPP} we know that each $v_n$ satisfies the sub-dynamic programming principle: for any $(s,x,i)\in[0,T]\times\R^d\times\I_m$ and $\rho\in\Tc^s$,
\beqs
v_n(s,x,i) &\leq & \sup_{\alpha\in\Ac_{s^+}}\inf_{u\in\Uc_s} \E\bigg[\int_s^{\rho'} f(X_t,I_t,u_t)dt + v_n(\rho',X_{\rho'},I_{\rho'}) \\
& & \hspace{2.5cm} - \; \sum_{n\in\N}c(X_{\tau_n'},I_{(\tau_n')^-},I_{\tau_n'})1_{\{s\leq\tau_n'<\rho'\}}\bigg].
\enqs
Since $v_n\leq v^-$, we get
\beq
\label{HalfDPP_v_n_v}
v_n(s,x,i) &\leq & \sup_{\alpha\in\Ac_{s^+}}\inf_{u\in\Uc_s} \E\bigg[\int_s^{\rho'} f(X_t,I_t,u_t)dt + v^-(\rho',X_{\rho'},I_{\rho'}) \\
& & \hspace{2.5cm} - \; \sum_{n\in\N}c(X_{\tau_n'},I_{(\tau_n')^-},I_{\tau_n'})1_{\{s\leq\tau_n'<\rho'\}}\bigg]. \notag
\enq
Letting $n\rightarrow\infty$ in \eqref{HalfDPP_v_n_v}, we finally obtain the sub-dynamic programming principle \eqref{DPP^-_v^-} for $v^-$. In a similar way we can prove \eqref{DPP^+_v^+}. Combining \eqref{DPP^-_v^-} and \eqref{DPP^+_v^+} with the equalities $v^-=V=v^+$, we end up with the dynamic programming principle for $V$.
%
\ep

\subsection{Elliott-Kalton formulation}
\label{Strategies}

We now describe the Elliott-Kalton formulation of the robust switching control problem, and we present in the next paragraph an example which shows that this is in general a different control problem than the robust feedback switching control problem studied here. As a by-product of this example, we will find a counterexample to uniqueness for equation \eqref{HJB}. Let us begin introducing the concept of non-anticipating strategy for the switcher. Firstly, we define a standard switching control, not necessarily of feedback form.

\begin{Definition}[Switching controls]
\label{D:switching_not_feedback}
Fix $s\in[0,T]$. We say that the double sequence $\alpha=(\tau_n,\iota_n)_{n\in\N}$ is a switching control starting at $s$ if:
\begin{itemize}
\item $\tau_n$ is an $\F^s$-stopping time, for any $n\in\N$, and
\[
s \ \leq \ \tau_0 \ \leq \ \cdots \ \leq \ \tau_n \ \leq \ \cdots \ \leq \ T.
\]
Moreover, $(\tau_n)_{n\in\N}$ satisfies the following property: for $\P$-a.e. $\omega\in\Omega$,
\[
\tau_n(\omega) \ = \ T, \qquad \text{for $n$ large enough}.
\]
\item $\iota_n\colon\Omega\rightarrow\I_m$ is $\Fc_{\tau_n}^s$-measurable, for any $n\in\N$.
\end{itemize}
$\overline\Ac_s$ denotes the family of all switching controls starting at $s$.
\end{Definition}

When using switching controls as defined above, the well-posedness of equation \eqref{SDE} becomes easier. In particular, we have the following result, whose standard proof is omitted.

\begin{Proposition}
\label{P:X_not_feedback}
Let Assumption {\bf (H1)} hold. For any $(s,x,i)\in[0,T]\times\R^d\times\I_m$, $\alpha\in\overline\Ac_s$, $u\in\Uc_s$, there exists a unique (up to indistinguishability) $\F^s$-adapted process $(X^{s,x,i;\alpha,u},I^{s,i;\alpha})=(X_t^{s,x,i;\alpha,u},I_t^{s,i;\alpha})_{s\leq t\leq T}$ to equation \eqref{SDE}. Moreover, estimate \eqref{EstimateX} holds.
\end{Proposition}

We can now introduce the concept of non-anticipating strategy for the switcher.

\begin{Definition}[Non-anticipating strategies]
Fix $s\in[0,T]$. We say that the map
\begin{align*}
\beta\colon\Uc_s&\longrightarrow\overline\Ac_s \\
u&\longmapsto\beta[u]=\big(\tau_n[u],\iota_n[u]\big)_{n\in\N}
\end{align*}
is a non-anticipating strategy starting at $s$ if
\[
\P\big[(\tau_n[u^1],\iota_n[u^1])1_{\{\tau_n[u^1]\leq t\}}=(\tau_n[u^2],\iota_n[u^2])1_{\{\tau_n[u^2]\leq t\}},\,\forall\,n\in\N\big] \ = \ 1
\]
whenever $\P(u_r^1=u_r^2,\,\forall\,r\in[s,t])=1$, for any $t\in[s,T]$ and $u^1,u^2\in\Uc_s$. $\Delta_s$ denotes the family of all non-anticipating strategies starting at $s$.
\end{Definition}

We can now define the corresponding value function:
\beqs
\hat V(s,x,i) &:=& \sup_{\beta\in\Delta_s}\inf_{u\in\Uc_s} J(s,x,i;\beta[u],u),
\enqs
for all $(s,x,i)\in[0,T]\times\R^d\times\I_m$. Notice that
\beq
\label{V<U}
V(s,x,i) &\leq & \hat V(s,x,i), \qquad \forall\,(s,x,i)\in[0,T]\times\R^d\times\I_m.
\enq
Under Assumptions {\bf (H1)} and {\bf (H2)}, we expect that $\hat V$ (as $V$) is a viscosity solution to equation \eqref{HJB}. Therefore, when {\bf (H3)} holds, by comparison, we have $V=\hat V$. However, if {\bf (H3)} is not assumed, the above inequality \eqref{V<U} might be strict at some $(s,x,i)\in[0,T]\times\R^d\times\I_m$. The following example illustrates this latter point.

\paragraph{Example.}
Fix $d=1$, $m=2$ so that $\I_2=\{1,2\}$, and take $U=\I_2$. Moreover, set $b(x,i,u)=-|i-u|$ and $\sigma\equiv0$. Notice that $b\in\{-1,0\}$. Since Assumption {\bf (H1)} is satisfied, from Proposition \ref{P:X_not_feedback} it follows that, for any $(s,x,i)\in[0,T]\times\R\times\I_2$, $\alpha\in\overline\Ac_s$, $u\in\Uc_s$, there exists a unique solution $(X^{s,x,i;\alpha,u},I^{s,i;\alpha})=(X_t^{s,x,i;\alpha,u},I_t^{s,i;\alpha})_{s\leq t\leq T}$ to equation \eqref{SDE}.

Set $g(x,i)=x$, $f\equiv0$, and $c\equiv0$. Our aim is now to determine the explicit form of $\hat V$ and $V$. To this end, it is convenient to give the following definition.

\begin{Definition}[Step controls]
\label{D:StepControls}
Fix $s\in[0,T]$. We say that $u$ is a step control starting at $s$ if there exists $n$ $\in$ $\N\backslash\{0\}$:
\begin{itemize}
\item $s=:t_0\leq\cdots\leq t_k\leq\cdots\leq t_n:=T$.
\item $\xi_k\colon\Omega\rightarrow U$ is $\Fc_{t_k}^s$-measurable, for any $k=0,\ldots,n-1$.
\end{itemize}
The control $u\colon[s,T]\times\Omega\rightarrow U$ is given by $u_t:=\sum_{k=0}^{n-1} \xi_k 1_{\{t_k\leq t<t_{k+1}\}}$. $\Uc_s^S$ denotes the family of all step controls starting at $s$.
\end{Definition}

Let us now determine the form of the function $\hat V$. Since the terminal payoff $g$ is strictly increasing and the drift $b$ is nonpositive, the aim of the switcher is to keep the system still. Having this in mind, we define, for every $\eps>0$, the strategy $\beta^\eps\in\Delta_s$, with $\beta^\eps[u]=(\tau_n^\eps[u],\iota_n^\eps[u])_{n\in\N}$ for all $u\in\Uc_s$, as follows.
\begin{itemize}
\item[(i)] For any $u_t=\sum_{k=0}^{n-1} \xi_k 1_{\{t_k\leq t<t_{k+1}\}}$ in $\Uc_s^S$, we set
\[
\big(\tau_k^\eps[u],\iota_k^\eps[u]\big) \ := \ (t_k,\xi_k), \qquad \forall\,k=0,\ldots,n-1.
\]
With this choice, $X_t^{s,x,i;\beta^\eps[u],u}=x$ for any $t\in[s,T]$ and $J(s,x,i;\beta^\eps[u],u)=x$.
\item[(ii)] For any $u\in\Uc_s\backslash\Uc_s^S$, it follows from the approximation result in \cite{krylov80}, Lemma 3.2.6, that there exists $u^\eps\in\Uc_s^S$ such that $\E[\int_s^T |u_t-u_t^\eps|dt]\leq\eps$. Then we define $\beta^\eps[u]:=\beta^\eps[u^\eps]$, where $\beta^\eps[u^\eps]$ has already been defined in item (i), since $u^\eps\in\Uc_s^S$. Therefore
\begin{align*}
&J(s,x,i;\beta^\eps[u],u) \ = \ \E\big[X_T^{s,x,i;\beta^\eps[u],u}\big] \ = \ x - \E\bigg[\int_s^T \big|I_T^{s,x,i;\beta^\eps[u],u} - u_t\big| dt\bigg] \\
&= \ x - \E\bigg[\int_s^T \big|I_T^{s,x,i;\beta^\eps[u^\eps],u} - u_t\big| dt\bigg] \ = \ x - \E\bigg[\int_s^T \big|u_t^\eps - u_t\big| dt\bigg] \ \geq \ x - \eps.
\end{align*}
\end{itemize}
In conclusion, we find, for every $\eps>0$,
\[
J(s,x,i;\beta^\eps[u],u) \ \geq \ x - \eps, \qquad \forall\,u\in\Uc_s,
\]
which implies $\inf_{u\in\Uc_s}J(s,x,i;\beta^\eps[u],u)\geq x-\eps$, and then $\hat V(s,x,i)\geq x-\eps$. From the arbitrariness of $\eps$, we obtain $\hat V(s,x,i)\geq x$. On the other hand, since $J(s,x,i;\beta[u],u)=\E[X_T^{s,x,i;\beta[u],u}]=x-\E[\int_s^T|I_t^{s,x,i;\beta[u],u}-u_t|dt]\leq x$, we deduce that
\[
\hat V(s,x,i) \ = \ g(x,i) \ = \ x, \qquad \forall\,(s,x,i)\in[0,T]\times\R\times\I_2.
\]
As a consequence of this result, we also have $\hat V(s,x,i)=\sup_{\beta\in\Delta_s}\inf_{u\in\Uc_s^S} J(s,x,i;\beta[u],u)$, for all $(s,x,i)\in[0,T]\times\R^d\times\I_m$.

Let us now find the expression of $V$. Fix $(s,x,i)\in[0,T]\times\R\times\I_2$ and $\alpha=(\tau_n,\iota_n)_{n\in\N}\in\Ac_s$. The aim of nature is to minimize the quantity $J(s,x,i;\alpha,u)$ over $\Uc_s$, which means to maximize the drift $b$, i.e., to keep it at the value $-1$. This can be done as follows. Define $u\in\Uc_s$, depending on $\alpha$, by
\[
u_t \ := \ (3-i)\,1_{\{s\leq t\leq\tau_0\}} + \sum_{n\in\N} (3-\iota_n)1_{\{\tau_n<t\leq\tau_{n+1}\}}, \qquad \forall\,t\in[s,T].
\]
Observe that, since $i,\iota_n\in\I_2$ then $3-i,3-\iota_n\in\I_2$; moreover, when $i=1$ then $3-i=2$, while if $i=2$ then $3-i=1$. Notice that, for $\P$-a.e. $\omega\in\Omega$ we have $I_t^{s,i;\alpha}(\omega)=3-u_t(\omega)$, for all $t\in[s,T]$ with $t\neq\tau_n(\omega)$, $n\in\N$. Therefore, $\P$-a.s.,
\[
b\big(X_t^{s,x,i;\alpha,u},I_t^{s,i;\alpha},u_t\big) \ = \ -\big|I_t^{s,i;\alpha} - u_t\big| \ = \ -1,
\]
for all $t\in[s,T]$, with $t\neq\tau_n$, $n\in\N$. It follows that, $\P$-a.s. we have $X_T^{s,x,i;\alpha,u}=x-(T-s)$. In other words, we obtain
\[
V(s,x,i) \ = \ x - (T - s), \qquad \forall\,(s,x,i)\in[0,T]\times\R\times\I_2.
\]
In conclusion, $V<\hat V$ on $[0,T)\times\R\times\I_2$. We finally observe that both $V$ and $\hat V$ are classical solutions to equation \eqref{HJB}, so that comparison does not hold. This is due to the fact that while Assumptions {\bf (H1)} and {\bf (H2)} hold, the no free loop property in {\bf (H3)} is not satisfied.

\begin{Remark}
{\rm
In the example above, because of the assumption that the switching costs are always zero ($c\equiv0$), it would be more natural, at least intuitively, to formulate the robust switching control problem as a classical two-player zero-sum stochastic differential game as in \cite{fleming_souganidis89}. In this latter setting, we recall from Theorem 2.6 in \cite{fleming_souganidis89} that the lower value function $V^{FS}$ (see Definition 1.4 in \cite{fleming_souganidis89}) is the unique viscosity solution to the lower Bellman-Isaacs equation:
\beq
\label{HJBI-}
\begin{cases}
-\dfrac{\partial w}{\partial t}(s,x) - \max_{i\in\I_2}\min_{u\in\I_2}\big[\Lc^{i,u}w(s,x)\big] \ = \ 0, \quad (s,x)\in[0,T)\times\R^d, \\
w(T,x) \ = \ x, \quad x\in\R^d,
\end{cases}
\enq
where $\Lc^{i,u}w(s,x)=-|i-u|D_xw(s,x)$. On the other hand, the upper value function $U^{FS}$ (see Definition 1.4 in \cite{fleming_souganidis89}) is the unique viscosity solution to the upper Bellman-Isaacs equation:
\beq
\label{HJBI+}
\begin{cases}
-\dfrac{\partial w}{\partial t}(s,x) - \min_{u\in\I_2}\max_{i\in\I_2}\big[\Lc^{i,u}w(s,x)\big] \ = \ 0, \quad (s,x)\in[0,T)\times\R^d, \\
w(T,x) \ = \ x, \quad x\in\R^d.
\end{cases}
\enq
By direct calculation, we see that $V$ satisfies \eqref{HJBI-}, so that it coincides with the lower value function $V^{FS}$ (this is expected from the results of \cite{fleming_souganidis89} and \cite{sirbu14b}, since $V$ is the sup/inf over feedback strategies/open-loop controls); while $\hat V$ satisfies \eqref{HJBI+}, therefore it coincides with the upper value function $U^{FS}$ (this is also not surprising, since $\hat V$ is the sup/inf over strategies/open-loop controls). Notice that in the present framework the Isaacs condition does not hold:
\[
\max_{i\in\I_2}\min_{u\in\I_2}[-|i-u|p] \ \neq \ \min_{u\in\I_2}\max_{i\in\I_2}[-|i-u|p], \qquad \forall\,p\in\R. 
\]
\ep
}
\end{Remark}

\vspace{9mm}

\small
\bibliographystyle{plain}
\bibliography{biblio}

\end{document}